\documentclass[12pt]{article}
\usepackage{amssymb,latexsym}
\usepackage{amsmath}
\usepackage{hyperref}
\usepackage{amsthm}
\usepackage{accents}
\usepackage{cite}
\usepackage{graphicx}
\usepackage{color}

\usepackage[normalem]{ulem}
\numberwithin{equation}{section}

\newtheorem{Thm}{Theorem}[section]
\newtheorem{Lem}[Thm]{Lemma}
\newtheorem{Def}[Thm]{Definition}

\newtheorem{Prop}[Thm]{Proposition}
\newtheorem{Rem}[Thm]{Remark}
\newtheorem{Cor}[Thm]{Corollary}
\newtheorem{Ex}[Thm]{Example}

\newcommand\de{\delta}
\newcommand\R{{\mathbb R}}
\newcommand\Sphere{{\mathbb S}}
\renewcommand\({\left(}
\renewcommand\){\right)}
\newcommand\hull{\operatorname{span}}
\newcommand{\photo}{P^{\:\!n}}
\renewcommand\d{\partial}
\newcommand\D{\nabla}
\newcommand\grad{\operatorname{grad}}
\newcommand\g{\gamma}
\renewcommand\a{\lambda}

\newcommand{\tr}{\operatorname{tr}}
\newcommand{\definedas}{\mathrel{\raise.095ex\hbox{\rm :}\mkern-5.2mu=}}
\newcommand{\asdefined}{\mathrel{=\mkern-5.2mu}\raise.095ex\hbox{\rm :}\;}
\newcommand{\surf}{\Sigma^{n-1}}
\newcommand{\schild}{Schwarz\-schild }
\newcommand{\slice}{M^{n}}
\newcommand{\Z}{\mathbb{Z}}
\newcommand{\Sphoto}{\overline{P}^{n}}  

\newcommand\beq{\begin{equation}}
\newcommand\eeq{\end{equation}}
\newcommand\ben{\begin{enumerate}}
\newcommand\een{\end{enumerate}}
\newcommand\bit{\begin{itemize}}
\newcommand\eit{\end{itemize}}
\newcommand*{\dt}[1]{ \accentset{\text{\large\bfseries .}}{#1}}
\newcommand{\Ric}{\operatorname{Ric}}
\newcommand{\Scal}{\operatorname{R}}

\title{Photon surfaces with equipotential time-slices}

\author{Carla Cederbaum\thanks{cederbaum@math.uni-tuebingen.de} \\Mathematics Department \\Eberhard Karls Universit\"at T\"ubingen\\ \\Gregory J. Galloway\thanks{galloway@math.miami.edu}
 \\Department of Mathematics \\ University of Miami}

\begin{document}

\date{}
\maketitle

\begin{abstract}
Photon surfaces are timelike, totally umbilic hypersurfaces of Lorentzian spacetimes. In the first part of this paper, we locally characterize all possible photon surfaces in a class of static, spherically symmetric spacetimes that includes (exterior) Schwarzschild, Reissner--Nordstr\"om, Schwarzschild-anti de Sitter, etc., in $n+1$ dimensions. In the second part, we prove that any static, vacuum, ``asymptotically isotropic'' $n+1$-dimensional spacetime that possesses what we call an ``equipotential'' and ``outward directed'' photon surface is isometric to the \schild spacetime of the same (necessarily positive) mass, using a uniqueness result by the first named author.
\end{abstract}

\section{Introduction}

One of the cornerstone results in the theory of black holes (in $3+1$ dimensions) is the static black hole uniqueness theorem, first due to Israel~\cite{Israel} for a single horizon and later due to Bunting and Masood-ul-Alam~\cite{BMuA} for multiple horizons, which establishes the uniqueness of the \schild spacetime among all static, asymptotically flat, black hole solutions to the vacuum Einstein equations. See Heusler's book~\cite{Heusler} and Robinson's review article~\cite{RobinsonReview} for a (then) complete list of references on further contributions, Simon's spinor proof recently described in Raulot's article~\cite[Appendix~A]{Raulot} and the recent article by Agostiniani and Mazzieri~\cite{AM} for newer approaches, both in the case of a single horizon.

A well-known, and intriguing, feature of (positive mass) \schild spacetime is the existence of a \emph{photon sphere}, namely, the timelike cylinder $P$ over the $\lbrace{r = (nm)^{\frac{1}{n-2}}, t=0\rbrace}$ $n-1$-sphere.  $P$ has the property of being \emph{null totally geodesic} in the sense that any null geodesic tangent to $P$ remains in $P$, i.e., $P$ traps all light rays tangent to it.

In~\cite{CDiss}, the first named author introduced, and made a study of the notion of a photon sphere for general static spacetimes (see also~\cite{CVE} in the spherically symmetric case).   Based on this study, in~\cite{CederPhoto}, she adapted Israel's argument (which requires the static lapse function to have nonzero gradient) to obtain a \emph{photon sphere uniqueness result}, thereby establishing the uniqueness of the \schild spacetime among all static, asymptotically flat solutions to the  vacuum Einstein equations which admit a single photon sphere.   Subsequent to that work, by adapting the argument of Bunting and Masood-ul-Alam~\cite{BMuA}, the authors~\cite{cedergal} were able to improve this result by, in particular, avoiding the gradient condition and allowing a priori multiple photon spheres.  For further results on photon spheres, in particular uniqueness results in the electro-vacuum case and the case of other matter fields, see for example~\cite{CedrGal2, GibbonsWarnick, yazadjiev, YazaLazov, YazaLazov2, Jahns, Shoom, Tomi, Tomi2, Yoshino}.

In this paper we will be concerned with the notion of \emph{photon surfaces} in spacetimes (see~\cite{CVE, Perlick} for 
slightly more general versions of this notion.) A photon surface in an $n+1$-dimensional spacetime $(\mathfrak{L}^{n+1}, \mathfrak{g})$ is a timelike hypersurface $\photo$ which is null totally geodesic, as described above.   As was shown in~\cite{CVE, Perlick}, a timelike hypersurface $\photo$ is a photon surface if and only if it is totally umbilic.
By definition, a \emph{photon sphere} in a \emph{static spacetime} is a photon surface $\photo$ along which the static lapse function $N$ is constant; see Section \ref{sec:prelim} for details.  While the \schild spacetime (in dimension $n+1$, $n \ge 3$, and with $r > (2m)^{\frac1{n-2}}$, $m>0$) admits a single photon sphere, it admits infinitely many photon surfaces of various types as briefly described in Section \ref{sec:spherical}.

In Section \ref{sec:spherical}, we derive the relevant ODE's describing spherically symmetric photon surfaces for a class of static, spherically symmetric spacetimes, which includes (exterior) Schwarzschild, Reissner--Nordstr\"om, and Schwarzschild-AdS. In the generic case, a photon surface in this setting is given by a formula for the derivative of the radius-time-profile curve $r = r(t)$; see Theorem \ref{thm:sphsymm}.  This formula is used in~\cite{CJM} to give a detailed qualitative description of all spherically symmetric photon surfaces in many (exterior) black hole spacetimes within class $\mathcal{S}$, including the (positive mass) \schild spacetime.  In addition, in Section \ref{sec:spherical}, we present a result which shows, for generic static, isotropic spacetimes which includes positive mass Schwarzschild and sub-extremal Reissner--Nordstr\"om, that, apart from some (partial) timelike hyperplanes, all photon surfaces are isotropic, see Theorem~\ref{thm:isotropic} and Corollary~\ref{coro:schwarzphoto}. As a consequence, one obtains a complete characterization of all photon surfaces in Reissner--Nordstr\"om and Schwarzschild.

Finally, in Section~\ref{sec:equi}, we obtain a new rigidity result pertaining to photon surfaces, rather than just to photon spheres. We prove that any static, vacuum, asymptotically isotropic spacetime possessing a (possibly disconnected) ``outward directed'' photon surface inner boundary with the property that the static lapse function $N$ is constant on each component of each time-slice $\surf(t)\definedas \photo\cap \lbrace{t = \text{const.}\rbrace}$ must necessarily be a \schild spacetime of positive mass, with the photon surface being one of the spherically symmetric photon surfaces in \schild classified in Section~\ref{sec:spherical}. We call such photon surfaces \emph{equipotential}. This generalizes static vacuum photon sphere uniqueness to certain photon surfaces and to higher dimensions.

The proof makes use of a new higher dimensional uniqueness result for the \schild spacetime  due to the first named author~\cite{ndimunique};  see Section \ref{sec:unique} for a statement.  This result generalizes in various directions the higher dimensional \schild uniqueness result of Gibbons et al.~\cite{GIZ}.  In particular, it does not a priori require the spacetime to be vacuum or static. A different proof of the result we use from~\cite{ndimunique} has since been given by Raulot~\cite{Raulot}, assuming that the manifolds under consideration are spin. These results rely on the rigidity case of a (low regularity version) of the Riemannian Positive Mass Theorem~\cite{SchoenYau,Witten,SYn,Miao,MSz,LeeLefloch}.

\section{Preliminaries}\label{sec:prelim}
The static, spherically symmetric $(n+1)$-dimensional \emph{\schild spacetime of mass $m\in\R$}, with $n\geq3$, is given by $(\overline{\mathfrak{L}}^{n+1}\definedas\R\times(\R^{n}\setminus \overline{B_{r_{m}}(0)}),\overline{\mathfrak{g}})$, where the Lorentzian metric $\overline{\mathfrak{g}}$ is given by
\begin{align}\label{schwarzmetric}
\overline{\mathfrak{g}}&=-\overline{N}^{2}dt^{2}+\overline{N}^{-2}dr^{2}+r^{2}\Omega,\quad 
\overline{N}=\left(1-\frac{2m}{r^{n-2}}\right)^{1/2},
\end{align}
with $\Omega$
denoting the standard metric on $\mathbb{S}^{n-1}$, and $r_{m}\definedas(2m)^{\frac{1}{n-2}}$ for $m>0$ and $r_{m}\definedas0$ for $m\leq0$, see also Tangherlini~\cite{Tangherlini}. For $m>0$, the timelike, cylindrical hypersurface $\Sphoto:=\lbrace r=(nm)^{\frac{1}{n-2}}\rbrace$ is called the \emph{photon sphere} of the \schild spacetime because any null geodesic (or ``photon'') $\gamma\colon\R\to\overline{\mathfrak{L}}^{n+1}$ that is tangent to $\Sphoto$ for some parameter $\tau_{0}\in\R$ is necessarily tangent to it for all parameters $\tau\in\R$. In particular, the \schild photon sphere is a timelike hypersurface ruled by null geodesics spiraling around the central black hole of mass $m>0$ ``at a fixed distance''. 

The \schild photon sphere can be seen as a special case of what is called a ``photon surface''~\cite{CVE,Perlick} in a general \emph{spacetime} (or smooth Lorentzian manifold): 

\begin{Def}[Photon surface]\label{def:photo-surf}
A timelike embedded hypersurface $\photo\hookrightarrow\mathfrak{L}^{n+1}$ in a spacetime $(\mathfrak{L}^{n+1},\mathfrak{g})$ is called a \emph{photon surface} if any null geodesic initially tangent to $\photo$ remains tangent to $\photo$ as long as it exists or in other words if $\photo$ is \emph{null totally geodesic}.
\end{Def}
The one-sheeted hyperboloids in the Minkowski spacetime (\schild spacetime with $m=0$) are also examples of photon surfaces, see Section \ref{sec:spherical}. It will be useful to know that, by an algebraic observation, being a null totally geodesic timelike hypersurface is equivalent to being an umbilic timelike hypersurface:
\begin{Prop}[\!\!{\cite[Theorem II.1]{CVE}, \cite[Proposition 1]{Perlick}}]\label{prop:umbilic}
Let $(\mathfrak{L}^{n+1},\mathfrak{g})$ be a spacetime and $\photo\hookrightarrow\mathfrak{L}^{n+1}$ an embedded timelike hypersurface. Then $\photo$ is a photon surface if and only if it is \emph{totally umbilic}, that is, if and only if its second fundamental form is pure trace.
\end{Prop}

As stated above, the \schild spacetime is ``static'', by the following definition.
\begin{Def}[Static spacetime]\label{def:static}
A spacetime $(\mathfrak{L}^{n+1},\mathfrak{g})$ is called \emph{(standard) static} if it is a warped product of the form
\begin{align}\label{static}
\mathfrak{L}^{n+1}&=\R\times \slice,\quad \mathfrak{g}=-N^{2}dt^{2}+g,
\end{align}
where $(\slice,g)$ is a smooth Riemannian manifold and $N\colon\slice\to\R^{+}$ is a smooth function called the \emph{(static) lapse function} of the spacetime.
\end{Def}

\begin{Rem}[Static spacetime cont., (canonical) time-slices]\label{rem:slice}
We will slightly abuse standard terminology and also call a spacetime static if it is a subset (with boundary) of a warped product static spacetime $(\R\times\slice,\mathfrak{g}=-N^{2}dt^{2}+g)$, $\mathfrak{L}^{n+1}\subseteq\R\times\slice$, to allow for inner boundary $\partial\mathfrak{L}$ not arising as a warped product. We will denote the \emph{(canonical) time-slices} $\lbrace{t=\text{const.}\rbrace}$ of a static spacetime $(\mathfrak{L}^{n+1},\mathfrak{g})$, $\mathfrak{L}^{n+1}\subseteq\R\times\slice$ by $\slice(t)$ and continue to denote the induced metric and (restricted) lapse function on $\slice(t)$ by $g$, $N$, respectively.
\end{Rem}

In the context of static spacetimes, we will use the following definition of ``photon spheres'', extending that of~\cite{CDiss,CederPhoto,CVE}. Consistently, the \schild photon sphere clearly is a photon surface in the \schild spacetime in this sense. 

\begin{Def}[Photon sphere]\label{def:photonsphere}
Let $(\mathfrak{L}^{n+1},\mathfrak{g})$ be a static spacetime, $\photo\hookrightarrow\mathfrak{L}^{n+1}$ a photon surface. Then $\photo$ is called a \emph{photon sphere} if the lapse function $N$ of the spacetime is constant along each connected component of $\photo$.
\end{Def}

For our discussions in Sections \ref{sec:spherical}, \ref{sec:equi}, we will make use of the following definitions.

\begin{Def}[Equipotential photon surface]\label{def:equipotential}
Let $(\mathfrak{L}^{n+1},\mathfrak{g})$ be a static spacetime, $\photo\hookrightarrow\mathfrak{L}^{n+1}$ a photon surface. Then $\photo$ is called \emph{equipotential} if the lapse function $N$ of the spacetime is constant along each connected component of each \emph{time-slice} $\surf\definedas\photo\cap\slice(t)$ of the photon surface.
\end{Def}

\begin{Def}[Outward directed photon surface]\label{def:outward}
Let $(\mathfrak{L}^{n+1},\mathfrak{g})$ be a static spacetime, $\photo\hookrightarrow\mathfrak{L}^{n+1}$ a photon surface arising as the inner boundary of $\mathfrak{L}^{n+1}$, $\photo=\partial\mathfrak{L}$, and let $\eta$ be the ``outward'' unit normal to $\photo$ (i.e. the normal pointing into $\mathfrak{L}^{n+1}$). Then $\photo$ is called \emph{ outward directed} if the $\eta$-derivative of the lapse function $N$ of the spacetime is positive, $\eta(N)>0$, along $\photo$.
\end{Def}

As usual, a spacetime $(\mathfrak{L}^{n+1},\mathfrak{g})$ is said to be \emph{vacuum} or to satisfy the \emph{Einstein vacuum equation} if
\begin{align}\label{EE}
\mathfrak{Ric}&=0
\end{align}
on $\mathfrak{L}^{n+1}$, where $\mathfrak{Ric}$ denotes the Ricci curvature tensor of $(\mathfrak{L}^{n+1},\mathfrak{g})$. For a static spacetime, the Einstein vacuum equation~\eqref{EE} is equivalent to the \emph{static vacuum equations}
\begin{align}\label{SMEvac1}
N\,{\Ric}&={\nabla}^2 N\\\label{SMEvac2}
\Scal&=0
\end{align}
on $\slice$, where $\Ric$, $\Scal$, and $\nabla^{2}$ denote the Ricci and scalar curvature, and the covariant Hessian of $(\slice,g)$, respectively. Combining the trace of~\eqref{SMEvac1} with~\eqref{SMEvac2}, one obtains the covariant Laplace equation on $\slice$,
\begin{align}\label{SMEvac3}
\triangle N&=0.
\end{align}
It is clear that, provided~\eqref{SMEvac1} holds,~\eqref{SMEvac2} and~\eqref{SMEvac3} can be interchanged without losing information. Of course, the \schild spacetime $(\R\times\overline{M}^{n},\overline{\mathfrak{g}})$ is vacuum and thus~\eqref{SMEvac1} and~\eqref{SMEvac3} hold for the \schild spatial metric $\overline{g}=\overline{N}^{-2}dr^{2}+r^{2}\Omega$  and lapse $\overline{N}$ on its canonical time-slice $\overline{M}^{n}=\R^{n}\setminus\overline{B_{r_{m}}(0)}$. 

Curvature quantities of a spacetime $(\mathfrak{L}^{n+1},\mathfrak{g})$ such as the Riemann curvature endomorphism $\mathfrak{Rm}$, the Ricci curvature tensor $\mathfrak{Ric}$, and the scalar curvature $\mathfrak{R}$ will be denoted in gothic print. The Lorentzian metric induced on a timelike embedded hypersurface $\photo\hookrightarrow\mathfrak{L}^{n+1}$ will be denoted by $p$, the (outward, see Definition \ref{def:outward}) unit normal by $\eta$, and the corresponding second fundamental form and mean curvature by $\mathfrak{h}$ and $\mathfrak{H}=\tr_{p}\mathfrak{h}$, respectively. With this notation, Proposition \ref{prop:umbilic} can be restated to state that a photon surface is characterized by
\begin{align}\label{eq:umbilic photo}
\mathfrak{h}&=\frac{\mathfrak{H}}{n}\,p.
\end{align}
To set sign conventions: $\mathfrak{h}(X,Y) = p({}^p\D_X\eta,Y)$ for vectors $X,Y$ tangent to $P$.

If the spacetime $(\mathfrak{L}^{n+1},\mathfrak{g})$ is static, its time-slices $\slice(t)$ have vanishing second fundamental form $K=0$ by the warped product structure, or, in other words, the time-slices are totally geodesic. The \emph{time-slices of a photon surface} $\photo\hookrightarrow\mathfrak{L}^{n+1}$ will be denoted by $\surf(t)\definedas\photo\cap\slice(t)$, with induced metric $\sigma=\sigma(t)$, second fundamental form $h=h(t)$, and mean curvature $H=H(t)=\tr_{\sigma(t)}h(t)$ with respect to the outward pointing unit normal $\nu =\nu(t)$. As an intersection of a totally geodesic time-slice and a totally umbilic photon surfaces, $\surf(t)$ is necessarily totally umbilic, and we have
\begin{align}\label{eq:umbilic surf}
h(t)&=\frac{H(t)}{n-1}\,\sigma(t).
\end{align}

Our choice of sign of the mean curvature is such that the mean curvature of $\mathbb{S}^{n-1}\hookrightarrow\R^n$ is positive  with respect to the outward unit normal in Euclidean space.

The following proposition will be useful to characterize photon surfaces in vacuum spacetimes. \begin{Prop}[{\!\!\cite[Proposition 3.3]{CederPhoto}}]
Let $n\geq2$ and let $(\mathfrak{L}^{n+1},\mathfrak{g})$ be a smooth semi-Riemannian manifold possessing a totally umbilic embedded hypersurface $\photo\hookrightarrow\mathfrak{L}^{n+1}$. If the semi-Riemannian manifold $(\mathfrak{L}^{n+1},\mathfrak{g})$ is Einstein, or in other words if $\mathfrak{Ric}=\Lambda \mathfrak{g}$ for some constant $\Lambda\in\R$, then each connected component of $\photo$ has constant mean curvature $\mathfrak{H}$ and constant scalar curvature
\begin{align}\label{eq:scalar}
\Scal_{p}&=(n+1-2\tau)\Lambda+\tau\,\frac{n-1}{n}\mathfrak{H}^{2},
\end{align}
where $\tau\definedas\mathfrak{g}(\eta,\eta)$ denotes the causal character of the unit normal $\eta$ to $\photo$. 
\end{Prop}

In particular, connected components of photon surfaces in vacuum spacetimes ($\Lambda=0$) have constant mean curvature and constant scalar curvature, related via
\begin{align}\label{eq:constant scalar}
\Scal_{p}&=\frac{n-1}{n}\,\mathfrak{H}^{2}.
\end{align}

We will now proceed to define and discuss the assumption of asymptotic flatness and asymptotic isotropy of static spacetimes.

\begin{Def}[Asymptotic flatness]\label{def:AF}
A smooth Riemannian manifold $(\slice,g)$ with $n\geq3$ is called \emph{asymptotically flat} if 
the manifold $\slice$ is diffeomorphic to the union of a (possibly empty) compact set and an open \emph{end} $E^{n}$ which is diffeomorphic to $\R^{n}\setminus \overline{B}$, $\Phi=(x^{i})\colon E^{n}\to\R^{n}\setminus \overline{B}$, where $B$ is some centered open ball in $\R^{n}$, and 
\begin{align}\label{AF}
(\Phi_{*}g)_{ij}-\delta_{ij}&=\mathcal{O}_{k}(r^{1-\frac{n}{2}-\varepsilon})\\\label{AFscal}
\Phi_{*}\!\Scal&=\mathcal{O}_{0}(r^{-n-\varepsilon})
\end{align}
for $i,j=1,\dots,n$ on $\R^{n}\setminus\overline{B}$ as $r\definedas \sqrt{(x^{1})^{2}+\dots+(x^{n})^{2}}\to\infty$ for some $k\in\Z$, $k\geq2$ and $\varepsilon>0$. Here, $\delta$ denotes the flat Euclidean metric, and $\delta_{ij}$ its components in the Cartesian coordinates $(x^{i})$.

A static spacetime $(\mathfrak{L}^{n+1}=\R\times \slice,\mathfrak{g}=-N^{2}dt^{2}+g)$ is called \emph{asymptotically flat} if its Riemannian base $(\slice,g)$ is asymptotically flat as a Riemannian manifold and, in addition, its lapse function satisfies
\begin{align}\label{NAF}
N-1&=\mathcal{O}_{k+1}(r^{1-\frac{n}{2}-\varepsilon})
\end{align}
on $\R^{n}\setminus\overline{B}$ as $r\to\infty$, with respect to the same coordinate chart $\Phi$ and numbers $k\in\Z$, $k\geq2$, $\varepsilon>0$.  We will abuse language and call $\mathfrak{L}^{n+1}\subseteq\R\times \slice$ \emph{asymptotically flat}, as long as $\mathfrak{L}^{n+1}$ has timelike inner boundary $\partial\mathfrak{L}$.
\end{Def}

One can expect a well-known result by Kennefick and \'O Murchadha~\cite{KM} to generalize to higher dimensions, which would assert that static vacuum asymptotically flat spacetimes are automatically ``asymptotically isotropic'' in suitable asymptotic coordinates. Here, we will resort to assuming asymptotic isotropy, leaving the higher dimensional generalization of this result to be dealt with elsewhere.

\begin{Def}[Asymptotic isotropy~\cite{ndimunique}]\label{def:AI}
A smooth Riemannian manifold $(\slice,g)$ of dimension $n\geq3$ is called \emph{asymptotically isotropic (of mass $m$)} if the manifold $\slice$ is diffeomorphic to the union of a (possibly empty) compact set and an open \emph{end} $E^{n}$ which is diffeomorphic to $\R^{n}\setminus \overline{B}$, $\Psi=(y^{i})\colon E^{n}\to\R^{n}\setminus \overline{B}$, where $B$ is some centered open ball in $\R^{n}$, and if there exists a constant $m\in\R$ such that
\begin{align}\label{eq:AI}
(\Psi_*g)_{ij}-(\widetilde{g}_{m})_{ij}&=\mathcal{O}_{2}(s^{1-n}),
\end{align}
for $i,j=1,\dots,n$ on $\R^{n}\setminus\overline{B}$ as $s\definedas \sqrt{(y^{1})^{2}+\dots+(y^{n})^{2}}\to\infty$, where 
\begin{align}\label{eq:schildiso1}
\widetilde{g}_{m}&\definedas\varphi^{\frac{4}{n-2}}_{m}(s)\,\delta,\\\label{eq:schildiso2}
\varphi_{m}(s)&\definedas1+\frac{m}{2s^{n-2}}
\end{align}
 denotes the spatial \schild metric in \emph{isotropic coordinates}. 

A static spacetime $(\mathfrak{L}^{n+1}=\R\times \slice,\mathfrak{g}=-N^{2}dt^{2}+g)$ is called \emph{asymptotically isotropic (of mass m)} if its Riemannian base $(\slice,g)$ is asymptotically isotropic of mass $m\in\R$ as a Riemannian manifold and, in addition, its lapse function $N$ satisfies
\begin{align}\label{eq:NAI}
N-\widetilde{N}_{m}&=\mathcal{O}_{2}(s^{1-n})
\end{align}
on $\R^{n}\setminus\overline{B}$ as $s\to\infty$, with respect to the same coordinate chart $\Psi$ and mass $m$. Here, $\widetilde{N}_{m}$ denotes the \schild lapse function in isotropic coordinates, given by
\begin{align}\label{eq:schildiso3}
\widetilde{N}_{m}(s)&\definedas \frac{1-\frac{m}{2s^{n-2}}}{1+\frac{m}{2s^{n-2}}}. 
\end{align}
As before, we will abuse language and call $\mathfrak{L}^{n+1}\subseteq\R\times \slice$ \emph{asymptotically isotropic} as long as it has timelike inner boundary.
\end{Def}

Here, we have rewritten the \schild spacetime, spatial metric, and lapse function in isotropic coordinates via the radial coordinate transformation
\begin{align}
r\asdefined s\,\varphi_{m}^{\frac{2}{n-2}}(s).
\end{align}
For $m>0$, this transformation bijectively maps $r\in(r_{m},\infty)\mapsto s\in(s_{m},\infty)$, with $s_{m}\definedas\left(\frac{m}{2}\right)^{\frac{1}{n-2}}$. For $m=0$, this transformation is the identity on $\R^{+}$, while for $m<0$, it only provides a coordinate transformation for $r$ suitably large, namely corresponding to $s>\left(\frac{\vert m\vert}{2}\right)^{\frac{1}{n-2}}$.

\begin{Rem}
A simple computation shows that the parameter $m$ in Definition \ref{def:AI} equals the ADM-mass of the Riemannian manifold $(\slice,g)$ defined in~\cite{ADM}.
\end{Rem}

\begin{Rem}
One can analogously define asymptotically flat and asymptotically isotropic Riemannian manifolds and static spacetimes with multiple ends $E^{n}_{l}$ and associated masses $m_{l}$. 
\end{Rem}

With these definitions at hand, let us point out that photon spheres are always  outward directed in static, vacuum, asymptotically isotropic spacetimes, a fact which is a straightforward generalization to higher dimensions of \cite[Lemma 2.6 and Equation (2.13)]{cedergal}:

\begin{Lem}\label{lem:outward}
Let $\photo\hookrightarrow\mathfrak{L}^{n+1}$ be a photon sphere in a static vacuum asymptotically flat spacetime $(\mathfrak{L}^{n+1},\mathfrak{g})$. Then $\photo$ is  outward directed.
\end{Lem}

\section{Photon surfaces in a class of static, spherically symmetric spacetimes}\label{sec:spherical}
In this section, we will give a local characterization of photon surfaces in a certain class $\mathcal{S}$ of static, spherically symmetric spacetimes $(\R\times \slice,\mathfrak{g})$, which  includes the $n+1$-dimensional (exterior) \schild spacetime. We will first locally characterize the spherically symmetric photon surfaces in $(\R\times \slice,\mathfrak{g})\in\mathcal{S}$ in Theorem \ref{thm:sphsymm} and then show  in Theorem~\ref{thm:isotropic} and in particular in Corollary~\ref{coro:schwarzphoto} that there are essentially no other photon surfaces in spacetimes $(\R\times \slice,\mathfrak{g})\in\mathcal{S}$. As mentioned in above, these results have been used in~\cite{CJM} to give a detailed description of all photon surfaces in many spacetimes in class $\mathcal{S}$, including the (positive mass) \schild spacetime.

The class $\mathcal{S}$ is defined as follows: Let $(\R\times \slice,\mathfrak{g})$ be a smooth Lorentzian spacetime such that
\begin{align}\label{def:ST}
\slice&= \mathcal{I}\times\Sphere^{n-1}\ni (r,\xi)
\end{align}
for an open interval $\mathcal{I}\subseteq(0,\infty)$, finite or infinite, and so that there exists a smooth, positive function $f\colon\mathcal{I}\to\R$ for which we can express the spacetime metric $\mathfrak{g}$ as
\begin{align}\label{def:STmetric}
\mathfrak{g}&=-f(r)dt^{2}+\frac{1}{f(r)}dr^{2}+r^{2}\Omega
\end{align}
in the global coordinates $t\in\R$, $(r,\xi)\in \mathcal{I}\times\Sphere^{n-1}$, where $\Omega$ denotes the canonical metric on~$\Sphere^{n-1}$ of area $\omega_{n-1}$. A Lorentzian spacetime $(\R\times \slice,\mathfrak{g})\in\mathcal{S}$ is clearly spherically symmetric and moreover naturally (standard) static via the hypersurface orthogonal, timelike Killing vector field $\d_t$.

\begin{Rem}
Note that we do not assume that spacetimes $(\R\times \slice,\mathfrak{g})\in\mathcal{S}$ satisfy any kind of Einstein equations or have any special type of asymptotic behavior towards the boundary of the radial interval $\mathcal{I}$, such as being asymptotically flat or asymptotically hyperbolic as $r\nearrow\sup\, \mathcal{I}$, or such as forming a regular minimal surface as $r\searrow \inf\, \mathcal{I}$.
\end{Rem}

\begin{Rem}\label{rem:tinv}
As $\d_{t}$ is a Killing vector field, the time-translation of any photon surface in a spacetime $(\R\times \slice,\mathfrak{g})\in\mathcal{S}$ will also be a photon surface in $(\R\times \slice,\mathfrak{g})$. As all spacetimes $(\R\times \slice,\mathfrak{g})\in\mathcal{S}$ are also time-reflection symmetric (i.e. $t \to -t$ is an isometry), the time-reflection of any photon surface in $(\R\times \slice,\mathfrak{g})$ will also be a photon surface in $(\R\times \slice,\mathfrak{g})$. 
 \end{Rem}

While the form of the metric~\eqref{def:STmetric} is certainly non-generic even among static, spherically symmetric spacetimes, the class $\mathcal{S}$ contains many important examples of spacetimes, such as the Minkowski and (exterior) \schild spacetime, the (exterior) Reissner--Nordstr\"om spacetime, the (exterior) Schwarzschild-anti de Sitter spacetime, etc.,
(in $n+1$ dimensions), each for a specific choice of $f$. 

Before we proceed with characterizing photon surfaces in spacetimes in this class $\mathcal{S}$, let us first make the following natural definition.

\begin{Def}
Let $(\R\times \slice,\mathfrak{g})\in\mathcal{S}$. A connected, timelike hypersurface $\photo\hookrightarrow(\R\times \slice,\mathfrak{g})$ will be called \emph{spherically symmetric} if, for each $t_{0}\in \R$ for which the intersection $\slice(t_0) := \photo\cap\lbrace t = t_{0}\rbrace\neq\emptyset$, there exists a radius $r_{0}\in \mathcal{I}$ (where $\slice=\mathcal{I}\times\Sphere^{n-1}$) such that
\begin{align}\label{def:sphsymmphoto}
\slice(t_0) &=\lbrace t_{0}\rbrace\times\lbrace r_{0}\rbrace\times\Sphere^{n-1}\subset \lbrace t_{0}\rbrace\times \slice.
\end{align}
A future timelike curve $\gamma\colon I\to \photo$, parametrized by arclength on some open interval $I\subset\R$, is called a \emph{radial profile} of $\photo$ if $\gamma'\in\hull\lbrace\d_{t},\d_{r}\rbrace\subset T_{\gamma'}(\R\times \slice)$ on $I$ and if the orbit of $\gamma$ under the rotation generates $\photo$.
\end{Def}

With this definition at hand, we will now prove the following lemma  which will  be used in the proof of Theorem \ref{thm:sphsymm}.

\begin{Lem}\label{lem:constant}
Let $(\R\times \slice,\mathfrak{g})\in\mathcal{S}$ and let $\photo\hookrightarrow(\R\times \slice,\mathfrak{g})$ be a spherically symmetric timelike hypersurface.  Assume $\photo\hookrightarrow(\R\times \slice,\mathfrak{g})$ has radial profile $\gamma\colon I\to \photo$, which may be written as, $\gamma(s)=(t(s),r(s),\xi_{*})\in\R\times\mathcal{I}\times\Sphere^{n-1}$ for some fixed $\xi_{*}\in\Sphere^{n-1}$.

If  $\photo\hookrightarrow(\R\times \slice,\mathfrak{g})$ is a photon surface, i.e. is totally umbilic with umbilicity factor $\a$, then the following first order  ODEs holds on $I$
\begin{align}
\label{eq:tdot}
\dot{t}  &= \frac{\a r}{f(r)}  \,,  \\
\label{eq:rdot}
(\dot{r})^2 &=  \a^2 r^2 - f(r)  \,,
\end{align}
where $\a$ is constant (and where $\dot{} = \frac{d}{ds}$).  Conversely, provided $\dot{r} \ne 0$,  
if the ODE's~\eqref{eq:tdot},~\eqref{eq:rdot} hold, with $\a$  constant, then $P$ is a photon surface with umbilicity factor~$\a$.
\proof To simplify notation we write $P$ for $\photo\hookrightarrow(\R\times \slice,\mathfrak{g})$ and $f$ for $f(r)$. 
As in Section \ref{sec:prelim}, let $p$ and $\mathfrak{h}$ denote the induced metric and second fundamental form of $P$, respectively. 
 
Set $e_0 = \dot{\g}$, and extend it to all of $P$ by making it invariant under the rotational symmetries.  Thus, $e_0$ is the future directed unit tangent vector field to $P$  orthogonal to each time-slice $\{t(s) = \text{const.}\}$, 
$s \in I$.  In terms of coordinates we have
\beq\label{eq:e0}
e_0 = \dot{t} \frac{\d}{\d t}  +   \dot{r} \frac{\d}{\d r} \,.
\eeq
Let $\eta$ be the outward pointing unit normal field to $P$.  From~\eqref{def:STmetric} and~\eqref{eq:e0} we obtain
\beq
\eta  =\frac{\dot{r}}{f}\d_{t}+\dot{t}f\d_{r} \,.
\eeq

\smallskip
\noindent
{\bf Claim:}  $P$ is a photon suface, with umbilicity factor $\a = \frac{f}{r}\dot{t}$, if and only if $e_0$ satisfies,
\beq\label{eq:covacc}
{}^p\D_{e_0}e_0 = \a \eta  \,.
\eeq

\proof[Proof of the claim]

Extend $e_0$ to an orthonormal basis $\{e_0, e_1, \dots, e_n\}$  in a neighborhood of an arbitrary point in $P$.  Thus, each $e_{I}$, $I = 1, ...,n$, where defined, is tangent to the time-slices.  A simple computation then gives,
\begin{align}
{}^p\D_{e_I}\eta &=\frac{\dot{r}}{f}\D_{e_I}\d_{t}+\dot{t}f\D_{e_I}\d_{r}  = \frac{\dot{t}f}{r} e_I  \,,
\end{align}
from which it follows that
\beq
\mathfrak{h}(e_I,e_J) = \a\, p({}^p\D_{e_I} \eta, e_J) = \a \de_{IJ} \,, \quad I,J = 1, ..., n \,,
\eeq
where $\de_{IJ}$ is the Kronecker delta, and
\beq\label{eq:alpha}
\a = \frac{f}{r}\dot{t}  \,.
\eeq
Similarly,
\beq
\mathfrak{h}(e_0,e_I) = \mathfrak{h}(e_I,e_0) = p({}^p\D_{e_I}\eta, e_0)) = 0 \,.
\eeq
Hence, 
\beq\label{eq:matrix}
\left[\mathfrak{h}(e_I,e_J)\right]_{I,J = 0,...,n}
\eeq
is a diagonal matrix with $\mathfrak{h}(e_I,e_I) = \a$ for $I = 1, ..., n$.   It remains to consider 
$\mathfrak{h}(e_0,e_0)$.

The profile curve $\g$, and its rotational translates, are `longitudes' in  the `surface of revolution' $P$.   As such, each is a unit speed geodesic in $P$, from which it follows that,
\beq
{}^p\D_{e_0}e_0 = \ell \eta
\eeq
for some scalar $\ell$.  This implies that,
\beq
\mathfrak{h}(e_0,e_0) = p({}^p\D_{e_0} \eta, e_0) = -p({}^p\D_{e_0}e_0, \eta)  = -\ell\,.
\eeq
From this and~\eqref{eq:matrix} we conclude that $P$ is a photon surface if and only if $\ell= \a = \frac{f}{r}\dot{t}$, which establishes the claim.

\smallskip
Using the coordinate expressions for $e_0$, $\eta$ and $\a$, a straight forward computation shows that~\eqref{eq:covacc}, with $\a = \frac{f}{r}\dot{t}$, is equivalent to
the following system of second order ODE's in the coordinate functions $t = t(s)$ and $r = r(s)$, 
\begin{align}
\ddot{t} + \frac{f'}{f} \dot{r} \dot{t} &= \frac{\dot{r}}{r} \dot{t} \label{eq:tddot}  \, ,\\
\ddot{r} + \frac{f f'}{2} \dot{t}^2 -\frac{f'}{2f}\dot{r}^2   &= \frac{(f\dot{t})^2}{r} \label{eq:rddot}  \,.
\end{align}
\end{Lem}

Now assume $P$ is a photon sphere with umbilicity factor $\a$, so that, in particular,~\eqref{eq:tddot} and~\eqref{eq:alpha} hold.   Treating~\eqref{eq:tddot} as a first order linear equation in $\dot{t}$, we have
\begin{align*}
\ddot{t} + \left(\frac{f'}{f} \dot{r}  -\frac{\dot{r}}{r}\right)\dot{t} = 0
\end{align*}
which, multiplying through by the integrating factor $\frac{f}{r}$, gives $\frac{d}{ds}\left(\frac{f}{r} \dot{t}\right) =  0$ so that~\eqref{eq:tdot} holds, with $\a = \frac{f}{r} \dot{t}>0$ a constant on $P$.   The assumption that $\g$ is parameterized with respect to arc length, gives
\beq\label{eq:unit}
f(\dot{t})^2  -  \ \frac{1}{f} (\dot{r})^2 = 1 \, .
\eeq
Together with~\eqref{eq:tdot}, we see that~\eqref{eq:rdot} also holds.  

Conversely, now assume that~\eqref{eq:tdot},~\eqref{eq:rdot} hold, with $\a = \frac{f}{r}\dot{t}=$ const., and, in addition, that $\dot{r} \ne 0$.  Differentiating~\eqref{eq:tdot} with respect to $s$, and then using~\eqref{eq:alpha} easily implies~\eqref{eq:tddot}.   Differentiating~\eqref{eq:rdot} with respect to $s$, then using~\eqref{eq:alpha} and dividing out by $\dot{r}$  gives,
\beq\label{eq:rddot2}
\ddot{r} + \frac{f'}{2} = \frac{(f\dot{t})^2}{r}  \, .
\eeq
Together with~\eqref{eq:unit} (which follows from~\eqref{eq:tdot} and~\eqref{eq:rdot}), this implies~\eqref{eq:rddot}.  We have shown that~\eqref{eq:tddot} and~\eqref{eq:rddot} hold, from which it follows that~\eqref{eq:covacc} holds with 
$\a = \frac{f}{r}\dot{t}$.   Invoking  the claim then completes the proof of Lemma \ref{lem:constant}.\qed

From Lemma \ref{lem:constant} we obtain the following.  

\begin{Thm}\label{thm:sphsymm}
Let $(\R\times \slice,\mathfrak{g})\in\mathcal{S}$ and let $\photo\hookrightarrow(\R\times \slice,\mathfrak{g})$ be a 
spherically symmetric timelike hypersurface. Assume that $\photo\hookrightarrow(\R\times \slice,\mathfrak{g})$ is a photon surface, with umbilicity factor $\a$, i.e.
\begin{align*}
\mathfrak{h}&=\a p
\end{align*}
where $p$ and $\mathfrak{h}$ are the induced metric and second fundamenatal form of 
$\photo\hookrightarrow(\R\times \slice,\mathfrak{g})$, 
respectively.   

Let $\gamma\colon I\to \photo$ be a radial profile for $\photo$ and write $\gamma(s)=(t(s),r(s),\xi_{*})\in\R\times\mathcal{I}\times\Sphere^{n-1}$ for some $\xi_{*}\in\Sphere^{n-1}$. Then $\a$ is a positive constant and \underline{\emph{either}} $r\equiv r_{*}$ along $\gamma$ for some $r_{*}\in \mathcal{I}$ at which the \emph{photon sphere condition} 
\begin{align}\label{eq:photonsphere}
f'(r_{*})r_{*}=2f(r_{*})
\end{align}
holds, $\a=\frac{\sqrt{f(r_{*})}}{r_{*}}$, and $(\photo,p)=(\R\times\Sphere^{n-1},-f(r_{*})dt^{2}+r^{2}_{*}\,\Omega)$ is a cylinder and thus a photon sphere, \underline{\emph{or}} $r=r(t)$ can globally be written as a smooth, non-constant function of $t$ in the range of $\gamma$ and $r=r(t)$ satisfies the \emph{photon surface ODE}
\begin{align}\label{eq:ODE3}
\left(\frac{dr}{dt}\right)^{2} &= \frac{f(r)^{2}\,(\a^2 r^2 - f(r))}{\a^{2} r^{2}}.
\end{align}

Conversely, whenever the photon sphere condition $f'(r_{*})r_{*}=2f(r_{*})$ holds for some $r_{*}\in \mathcal{I}$, then the cylinder $(\photo,p)=(\R\times\Sphere^{n-1},-f(r_{*})dt^{2}+r^{2}_{*}\,\Omega)$ is a photon sphere in $(\R\times \slice,\mathfrak{g})$ with umbilicity factor $\a=\frac{\sqrt{f(r_{*})}}{r_{*}}$. Also, any smooth, non-constant solution $r=r(t)$ of~\eqref{eq:ODE3} for some constant $\a>0$ gives rise to a photon surface in $(\R\times \slice,\mathfrak{g})$ with umbilicity factor $\a$.  
\end{Thm}

\proof
From Lemma \ref{lem:constant}, we know that $\a$ is a positive constant. Moreover, we know that 
$t = t(s)$ and $r = r(s)$ satisfy equations~\eqref{eq:tdot} and~\eqref{eq:rdot}.

In the case when $r\equiv r_{*}$ for some constant $r_*$, these equations immediately imply,
\beq
\dot{t} = \frac{1}{\sqrt{f(r_{*})}}, \quad \a=\frac{\sqrt{f(r_{*})}}{r_{*}} \,.
\eeq
Further,~\eqref{eq:rddot}  implies
\beq
f'(r_{*})r_{*} = 2f(r_{*})  \,.
\eeq

In the general case, Equations~\eqref{eq:tdot} and~\eqref{eq:rdot} clearly imply~\eqref{eq:ODE3}.   The converse statements are easily obtained from~\eqref{eq:ODE3} and the unit speed condition~\eqref{eq:unit}.\qed

\begin{Rem}\label{ref:tinv2}
In view of Remark \ref{rem:tinv}, note that in the `either' case, the photon sphere is time-translation and time-reflection invariant in itself. In the `or' case, note that the photon surface ODE~\eqref{eq:ODE3} is time-translation and time-reflection invariant and will thus allow for time-translated and time-reflected solutions corresponding to the same $\a>0$.
\end{Rem}

\begin{Ex}\label{ex:Minkowski}
Choosing $(\R\times \slice,\mathfrak{g})=(\R^{1,n},\mathfrak{m})$, where $\mathfrak{m}$ is the Minkowski metric and $f\colon(0,\infty)\to\R\colon r\mapsto1$, the photon sphere condition cannot be satisfied for any $r_{*}\in(0,\infty)$ so that every spherically symmetric photon surface in the Minkowski spacetime must satisfy the ODE~\eqref{eq:ODE3} which reduces to
\begin{align*}
\left(\frac{dr}{dt}\right)^{2} &= \frac{\a^2 r^2 - 1}{\a^{2} r^{2}}\quad\Leftrightarrow\quad r(t)=\sqrt{\a^{-2}+(t-t_{0})^{2}}\text{ for some }t_{0}\in\R
\end{align*}
and describes the rotational one-sheeted hyperboloids of radii $\a^{-1}$ for any $0<\a<\infty$.
\newpage
This is of course consistent with the well-known fact that the only timelike, totally umbilic hypersurfaces in the Minkowski spacetime are, apart from (parts of) timelike hyperplanes, precisely (parts of) these hyperboloids and their spatial translates, the formula for which explicitly displays the time-translation and time-reflection invariance of the photon surface characterization problem.
\end{Ex}

\begin{figure}[h]
\begin{center}
\mbox{
\includegraphics[width=5.85in]{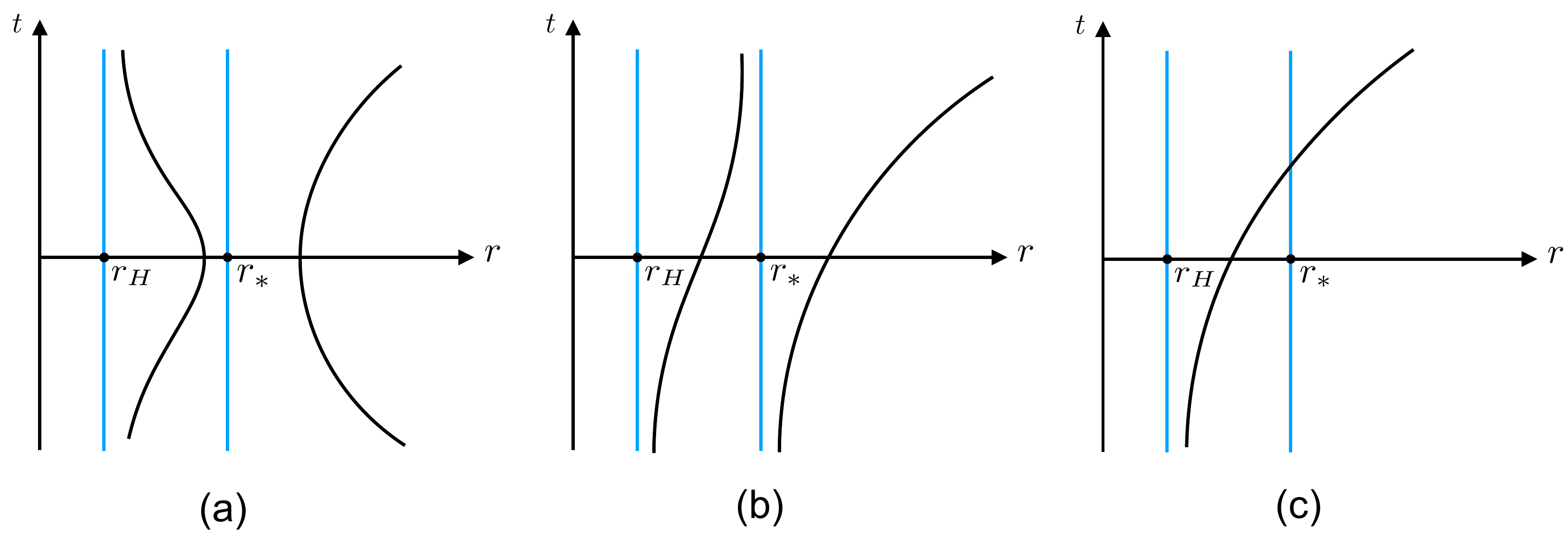}
}
\end{center}
\vspace{-1ex}\caption{Profile curves for all types of spherically symmetric photon surfaces in \schild spacetime grouped according to umbilicity factor; see~\cite{CJM} for details.}
\label{Photon}
\end{figure}

Note that the photon sphere condition is satisfied precisely at the well-known photon sphere radius $r_{*}=(nm)^{\frac{1}{n-2}}$ in the $n+1$-dimensional \schild spacetime where $f(r)=1-\frac{2m}{r^{n-2}}$ for $m>0$ and $r> r_H =(2m)^{\frac{1}{n-2}}$ and there is no photon sphere radius for $m\leq0$ and $r>0$.   
While there are no other photon spheres, there are many non-cylindrical photon surfaces in the \schild spacetime.  The analysis in~\cite{CJM}, based on Theorem~\ref{thm:sphsymm}, shows that, up to time translation and time reflection (cf. Remarks \ref{rem:tinv} and \ref{ref:tinv2}), there are five classes of non-cylindrical, spherically symmetric photon surfaces in the (exterior) positive mass \schild spacetime (as well as in many other (exterior) black hole spacetimes in class $\mathcal{S}$); the profile curves for representatives from each class are depicted in Figure \ref{Photon}. 

\vspace{-1ex}
In each case, they approach asymptotically the event horizon $r = r_H$ and/or the photon sphere $r = r_*$ and/or become asymptotically null at infinity in $(t,r)$-coordinates. An analysis of the behavior of the asymptotics of the non-cylindrical, spherically symmetric photon is performed for both Schwarzschild and many other (exterior) black hole spacetimes in class $\mathcal{S}$ in (generalized) Kruskal--Szekeres coordinates in~\cite{CW}. There, it is found that the photon surfaces appearing to approach $r = r_H$ in $(t,r)$-coordinates in fact cross the event horizon, while those approaching the photon sphere $r = r_*$ or asymptotically become null in $(t,r)$-coordinates do so in (generalized) Kruskal--Szekeres coordinates, too.

Using quite different methods, in~\cite{Perlicketal}, the same types of photon surfaces are found in a $2+1$-dimensional spacetime obtained by dropping an angle coordinate from $3+1$-dimensional Schwarzschild of positive mass.

The question naturally arises: What about non-spherically symmetric photon surfaces?  This is addressed in the following theorem, see in particular Corollary~\ref{coro:schwarzphoto}.
\begin{Thm}\label{thm:isotropic}  
Let $n\geq3$, $I\subseteq\R^{+}$ an open interval, $D^{n}\definedas \lbrace{y\in\R^{n}\,\vert\,\vert y\vert=s\in I\rbrace}$, and let $\widetilde{N},\psi\colon I\to\R^{+}$ be smooth, positive functions. Consider the \emph{static, isotropic spacetime}
\begin{align}
\left(\R\times D^{n},\mathfrak{g}=-\widetilde{N}^{2}dt^{2}+\psi^{2}\,\delta\right)
\end{align}
of lapse $\widetilde{N}=\widetilde{N}(s)$ and conformal factor~$\psi=\psi(s)$. We write $\widetilde{g}\definedas \psi^{2}\,\delta$. A timelike hypersurface $\photo$ in $(\R\times D^{n},\mathfrak{g})$ is called \emph{isotropic} if $\photo\cap \lbrace{t=\text{const.}\rbrace}=\mathbb{S}^{n-1}_{s(t)}(0)\subset D^{n}$ for some radius $s(t)\in I$ for every $t$ for which $\photo\cap \lbrace{t=\text{const.}\rbrace}\neq\emptyset$. A \emph{(partial) centered vertical hyperplane} in $(\R\times D^{n},\mathfrak{g})$ is the restriction of a timelike hyperplane in the Minkowski spacetime containing the $t$-axis to $\R\times D^{n}$, i.e. a set of the form
\begin{align}
\lbrace{(t,y)\in\R\times D^{n}\,\vert\, y\cdot u=0\rbrace}
\end{align}
for some fixed Euclidean unit vector $u\in\R^{n}$, where $\cdot$ denotes the Euclidean inner product. Centered vertical hyperplanes are totally geodesic in $(\R\times D^{n},\mathfrak{g})$. 

Assume furthermore that the functions $\widetilde{N}$ and $\psi$ satisfy
\begin{align}\label{eq:unless}
\frac{\widetilde{N}'(s)}{\widetilde{N}(s)}&\neq\frac{\psi'(s)}{\psi(s)}
\end{align}
for all $s\in I$. Then any photon surface in $(\R\times D^{n},\mathfrak{g})$ is either (part of) an isotropic photon surface or (part of) a centered vertical hyperplane.
\end{Thm}

\begin{Cor}\label{coro:schwarzphoto} 
Let $n\geq3$, $m>0$,  and consider the $n+1$-dimensional \schild spacetime of mass $m$. Then any connected photon surface is either (part of) a centered vertical hyperplane as described above or (part of) a spherically symmetric photon surface as described in Theorem~\ref{thm:sphsymm}.
\end{Cor}
\proof \emph{(of Corollary~\ref{coro:schwarzphoto})} Recall the isotropic form of the Schwarzschild spacetime~\eqref{eq:schildiso1},~\eqref{eq:schildiso2}, and~\eqref{eq:schildiso3}, with $I=(s_{m},\infty)$, and note that~\eqref{eq:unless} corresponds to $s^{n-2}\neq\frac{m}{2(n-1)}$ (which can be quickly seen when exploiting $\widetilde{N}=\frac{2-\varphi}{\varphi}$, $\psi=\varphi^{\frac{2}{n-2}}$). This, however, is automatic as $(\frac{m}{2(n-1)})^{\frac{1}{n-2}}<s_{m}=(\frac{m}{2})^{\frac{1}{n-2}}<s$.

\begin{Rem}\label{rem:nowhereconfflat}
Condition~\eqref{eq:unless} can be interpreted geometrically as follows: If 
$\frac{\widetilde{N}'(s)}{\widetilde{N}(s)}=\frac{\psi'(s)}{\psi(s)}$
holds for $s\in J$, $J\subseteq I$ an open interval, then there exists a positive constant $A>0$ such that $\widetilde{N}(s)=A\psi(s)$ for all $s\in J$, where we used that $\widetilde{N},\psi>0$. This shows that $\mathfrak{g}=-\widetilde{N}^{2}dt^{2}+\psi^{2}\delta=\psi^{2}(-A^{2}dt^{2}+\delta)$ or in other words the static, isotropic spacetime $(\R\times(J\times\mathbb{S}^{n-1}),\mathfrak{g})$ is globally conformally flat and hence possesses additional  photon surfaces corresponding to the totally geodesic timelike hyperplanes that do not contain the $t$-axis and to the spatially translated totally umbilic rotational one-sheeted hyperboloids of the (time-rescaled) Minkowski spacetime, see  Example~\ref{ex:Minkowski}.

On the other hand, reconsidering the proof of Theorem~\ref{thm:isotropic}, one finds that the assumption~\eqref{eq:unless} is not (fully) necessary (it is only needed for the reasoning after~\eqref{eq:mixedh} and~\eqref{eq:mixedhspherical} for the planar and the spherical case, respectively). Hence, Theorem~\eqref{thm:isotropic} gives a full characterization of photon surfaces in nowhere locally conformally flat static, isotropic spacetimes.
\end{Rem}

\begin{Rem}\label{rem:consistent}
A static, isotropic spacetime $(\R\times D^{n},-\widetilde{N}^{2}dt^{2}+\psi^{2}\,\delta)$ can be globally rewritten as a spacetime of class $\mathcal{S}$ if and only if $\widetilde{N}^{2}(s)=(1+\frac{s\psi'(s)}{\psi(s)})^{2}>0$ for all $s\in I$ by setting $r(s)\definedas s\psi(s)$ and $f(r)\definedas \widetilde{N}(s(r))$, where $s=s(r)$ denotes the inverse function of $r=r(s)$. In this case, the photon sphere and photon surface conditions on the isotropic radius profile $s=S_{*}$ and $s=S(t)$, Equations~\eqref{eq:photonsphereisotropic} and~\eqref{eq:photonsurfaceisotropic}, reduce to the much simpler photon sphere and photon surface conditions for the area radius profile $r=r_{*}$ and $r=r(t)$, Equations \eqref{eq:photonsphere} and \eqref{eq:ODE3}, respectively. 

Conversely, a spacetime of class $\mathcal{S}$ can always be locally rewritten in isotropic form by picking a suitable $r_0\in\mathcal{I}$ (or $r_0\in\overline{\mathcal{I}}$) and setting 
\begin{align*}
s=s(r)\definedas\exp\left(\int^r_{r_0} \left(\rho\sqrt{(f(\rho)}\right)^{-1}\,d\rho\right),
\end{align*}
 $\psi(s)\definedas\frac{r(s)}{s}$ and $\widetilde{N}(s)\definedas\sqrt{f(r(s))}$, where $r=r(s)$ denotes the inverse of $s=s(r)$.

The main reason for switching into the isotropic picture lies in the spatial conformal flatness allowing to easily describe centered vertical hyperplanes and to exclude photon surfaces that are not centered vertical hyperplanes nor isotropic.
\end{Rem}

\proof[Proof of Theorem \ref{thm:isotropic}] Let $\photo$ be a connected photon surface in a static, isotropic spacetime $(\R\times D^{n},\mathfrak{g}=-\widetilde{N}^{2}dt^{2}+\psi^{2}\delta)$. As before, set $\slice(t)\definedas\lbrace{t=\text{const.}\rbrace}$. Let $T\definedas\lbrace{t\in\R\,\vert\,\photo\cap \slice(t)\neq\emptyset\rbrace}$ and note that $T$ is an open, possibly infinite, interval. Set $\surf(t)\definedas \photo\cap \slice(t)$ for $t\in T$. As timelike and spacelike submanifolds are always transversal, $\surf(t)$ is a smooth surface. Furthermore, $\surf(t)$ is umbilic in $\slice(t)$ by time-symmetry of $\slice(t)$, or in other words because the second fundamental form of $\slice(t)$ in a static spacetime vanishes. As $(\slice(t),\widetilde{g})$ is conformally flat, exploiting the conformal invariance of umbilicity, the only umbilic hypersurfaces in  $(\slice(t),\widetilde{g})$ are the conformal images of pieces of Euclidean round spheres and pieces of Euclidean hyperplanes. Slightly abusing notation and denoting points in the spacetime by their isotropic coordinates, by continuity and by connectedness of $\photo$, $\left(\surf(t)\right)_{t\in T}$ is thus either a family of pieces of spheres 
\begin{align}\label{eq:spheres}
\surf(t)\subseteq\lbrace{ y\in\R^{n}\,\vert\,\vert y-c(t)\vert = S(t) \rbrace}
\end{align} 
with centers $c(t)\in\R^{n}$ and radii $S(t)>0$ for all $t\in T$, or a family of pieces of hyperplanes 
\begin{align}\label{eq:hyperplanes}
\surf(t)\subseteq\lbrace{y\in\R^{n}\,\vert\,y\cdot u(t) = a(t) \rbrace}
\end{align}
for some $\delta$-unit normal vectors $u(t)\in\R^{n}$ and altitudes $a(t)\in\R$ for all $t\in T$, where $\vert\cdot\vert$ and $\cdot$ denote the Euclidean norm and inner product, respectively. The (outward, where appropriate) unit normal $\eta$ to $\photo$ can be written as $\eta=\alpha\nu+\beta\partial_{t}$, with $\alpha>0$, recalling that $\nu$ denotes the (outward, where appropriate) unit normal to $\surf(t)$ in $\slice(t)$. By time-symmetry of $\slice(t)$, the second fundamental form $\mathfrak{h}$ of $\photo$ in the spacetime restricted to $\surf(t)$ can be expressed in terms of the second fundamental form $h$ of $\surf(t)$ in $\slice(t)$ via
$\mathfrak{h}\vert_{T\surf(t)\times T\surf(t)}=\alpha h$. By umbilicity, $\mathfrak{h}=\lambda p$, $p$ denoting the induced metric on $\photo$, this implies 
\begin{align}\label{eq:humbilic}
h&=\frac{\lambda}{\alpha}\,\sigma,
\end{align} 
where $\sigma$ is the induced metric on $\surf(t)$. We will treat the planar and the spherical cases separately. We will denote $t$-derivatives by $\dot{}$ and $s$-derivatives by $'$. 

\paragraph*{Planar case:} Let $\surf(t)$ be as in~\eqref{eq:hyperplanes} for all $t\in T$. We will show that $a(t)=0$ and $\dot{u}(t)=0$ for all $t\in T$, and moreover that $\lambda=0$ along $\photo$. This then implies that $\photo$ is contained in a centered vertical hyperplane with unit normal $\eta=\nu=\psi^{-1}(s)u^{i}\partial_{y^{i}}$ and moreover that all centered vertical hyperplanes are totally geodesic as every centered vertical hyperplane can be written in this form. 

For each $t\in T$, extend $u(t)$ to a $\delta$-orthonormal basis $\lbrace{e_{1}(t)=u(t),e_{2}(t),\dots e_{n}(t)\rbrace}$ of $\R^{n}$ in such a way that $e_{I}(t)$ is smooth in $t$ for all $I=2,\dots,n$. Then clearly $X_{I}(t,y)\definedas e_{I}^{k}(t)\partial_{y^{k}}$ is tangent to $\photo$ for all $I=2,\dots n$ and $\lbrace{X_{I}(t,\cdot)\rbrace}_{I=2}^{n}$ is an orthogonal frame for $\surf(t)$ with respect to $\widetilde{g}$ by conformal flatness. To find the missing (spacetime-)orthogonal tangent vector to $\photo$, consider a curve $\mu(t)=(t,y(t))$ in $\photo$ with tangent vector $\dot{\mu}(t)=\partial_{t}+\dot{y}^{i}(t)\partial_{y^{i}}$. Let capital latin indices run from $2,\dots,n$.  Now decompose $\dot{y}(t)=\rho(t)u(t)+\xi^{I}(t)e_{I}(t)\in\R^{n}$. By~\eqref{eq:hyperplanes}, we find $\rho(t)=\dot{y}(t)\cdot u(t)=\dot{a}(t)-y(t)\cdot\dot{u}(t)$. Hence, a future pointing tangent vector to $\photo$ orthogonal in the spacetime to all $X_{I}$ is given by
\begin{align}\label{eq:X1}
X_{1}(t,y)&\definedas \partial_{t}+\left(\dot{a}(t)-y\cdot\dot{u}(t)\right)u^{i}(t)\partial_{y^{i}},
\end{align}
so that we have constructed a smooth orthogonal tangent frame $\lbrace{X_{i}\rbrace}_{i=1}^{n}$ for $\photo$. Hence we can compute the (spacetime) unit normal to $\photo$ to be
\begin{align}\label{eq:etaplanar}
\eta(t,y)&=\frac{\frac{\psi(s)}{\widetilde{N}(s)}\left(\dot{a}(t)-y\cdot\dot{u}(t)\right)\partial_{t}+\frac{\widetilde{N}(s)}{\psi(s)}\,u^{i}(t)\partial_{y^{i}}}{\sqrt{\widetilde{N}^{2}(s)-\psi^{2}(s)\left(\dot{a}(t)-y\cdot\dot{u}(t)\right)^{2}}}.
\end{align}
In other words, using that $\nu(t,y)=\psi^{-1}(s)\,u^{i}(t)\partial_{y^{i}}$, we have
\begin{align}\label{eq:alphaeta}
\alpha(t,y)&=\frac{\widetilde{N}(s)}{\sqrt{\widetilde{N}^{2}(s)-\psi^{2}(s)\left(\dot{a}(t)-y\cdot\dot{u}(t)\right)^{2}}},\\\label{eq:betaeta}
\beta(t,y)&=\frac{\psi(s)\left(\dot{a}(t)-y\cdot\dot{u}(t)\right)}{\widetilde{N}(s)\sqrt{\widetilde{N}^{2}(s)-\psi^{2}(s)\left(\dot{a}(t)-y\cdot\dot{u}(t)\right)^{2}}}.
\end{align}
We are now in a position to compute the second fundamental forms explicitly and take advantage of the umbilicity of $\photo$. Using $u(t)\cdot \dot{e}_{J}(t)=-\dot{u}(t)\cdot e_{J}(t)$ for all $t\in T$, the condition $\mathfrak{h}(X_{1},X_{J})=0$ gives
\begin{align}
-\dot{u}(t)\cdot e_{J}(t)+\frac{1}{s}\left\{\frac{\psi'(s)}{\psi(s)}-\frac{\widetilde{N}'(s)}{\widetilde{N}(s)}\right\}\left(\dot{a}(t)-y\cdot\dot{u}(t)\right)y\cdot e_{J}(t)&=0
\end{align}
for $J=2,\dots,n$. As $\lbrace{u(t),e_{J}(t)\rbrace}_{J=2}^{n}$ is a $\delta$-orthonormal frame, this is equivalent to
\begin{align}\label{eq:mixedh}
\frac{1}{s}\left\{\frac{\psi'(s)}{\psi(s)}-\frac{\widetilde{N}'(s)}{\widetilde{N}(s)}\right\}\left(\dot{a}(t)-y\cdot\dot{u}(t)\right)\left(y-a(t)\,u(t)\right)&=\dot{u}(t).
\end{align}
As $\surf(t)$ has dimension $n-1$,~\eqref{eq:mixedh} tells us that $\dot{u}(t)=0$ for all $t\in T$ by linear dependence considerations (or otherwise if the term in braces $\{\dots\}$ vanishes). Hence,~\eqref{eq:mixedh} simplifies to
\begin{align}
\left\{\frac{\psi'(s)}{\psi(s)}-\frac{\widetilde{N}'(s)}{\widetilde{N}(s)}\right\}\,\dot{a}(t)\,\left(y-\,a(t)\,u\right)&=0
\end{align}
so that, for a given $t\in T$, again using that $\surf(t)$ has dimension $n-1$  and linear dependence considerations, we find $\dot{a}(t)=0$ if the term in braces $\{\dots\}$ does not vanish along $\photo$, i.e. when assuming~\eqref{eq:unless}. 

Let us now compute the umbilicity factor $\lambda$, exploiting that $u$ and $a$ are constant. Note that $e_{I}$ is also constant, and~\eqref{eq:alphaeta} and~\eqref{eq:betaeta} reduce to $\alpha=1$ and $\beta=0$, and, using~\eqref{eq:etaplanar} and~\eqref{eq:X1}, we get $\eta=\nu$ and $X_{1}=\partial_{t}$. As $e_{I}$ is independent of $y$, we find 
\begin{align}
h(X_{I},X_{J})&=\frac{\psi'(s)}{s}\,a\,\delta_{IJ}
\end{align}so that by~\eqref{eq:humbilic}, the photon surface umbilicity factor $\lambda$ satisfies
\begin{align}\label{eq:lambdah}
\lambda(t,y) &=\lambda(y)=\frac{\psi'(s)}{s\psi^{2}(s)}\,a
\end{align}
and is in particular independent of $t$. From $\mathfrak{h}(X_{1},X_{1})=\lambda p(X_{1},X_{1})$, we find
\begin{align}\label{eq:lengthy}
\begin{split}
\lambda(y) &=\frac{\widetilde{N}'(s)}{s\widetilde{N}(s)\psi(s)}\, a.
\end{split}
\end{align}
Thus,~\eqref{eq:lambdah} and~\eqref{eq:lengthy} combine to
\begin{align}
\lambda(y)
&=\frac{\psi'(s)}{s\psi^{2}(s)}\,a=\frac{\widetilde{N}'(s)}{s\widetilde{N}(s)\psi(s)}\,a
\end{align}
which implies $a=0$ and indeed $\lambda(y)=\lambda=0$ is also independent of $y$, when assuming that~\eqref{eq:unless} holds along $\photo$. This shows that centered vertical hyperplanes are totally geodesic and that any photon surface $\photo$ as in~\eqref{eq:hyperplanes} along which~\eqref{eq:unless} holds is (part of) a centered vertical hyperplane.

\paragraph*{Spherical case:} Let $\surf(t)$ be as in~\eqref{eq:spheres} for all $t\in T$. We will show that $c(t)=0$ for all $t\in T$. This then implies that $\photo$ is contained in an isotropic photon surface as desired, namely in a photon sphere with isotropic radius $s=S_{*}$ satisfying~\eqref{eq:photonsphereisotropic} or with isotropic radius profile $s=S(t)$ as in~\eqref{eq:photonsurfaceisotropic}. We will use the abbreviation
\begin{align}\label{eq:defu}
u(t,y)\definedas\frac{y-c(t)}{S(t)}
\end{align}
to reduce notational complexity.

A straightforward computation shows that the outward unit normal $\nu$ to $\surf(t)$ in $\slice(t)$ is given by 
\begin{align}
\nu&=\psi^{-1}(s)\,u^{i}(t,y)\,\partial_{y^{i}}.
\end{align}
Now choose a smooth $\delta$-orthonormal system of vectors $e_{I}(t,y)$ locally along $\photo$ such that $e_{I}(t,y)\cdot u(t,y)=0$ for all $(t,y)\in\photo$ and set $X_{I}(t,y)\definedas e_{I}^{k}(t,y)\partial_{y^{k}}$ for all $(t,y)\in\photo$ and all $I=2,\dots n$ so that $\lbrace{X_{I}(t,\cdot)\rbrace}_{I=2}^{n}$ is an orthogonal frame for $\surf(t)$ with respect to $\widetilde{g}$ by conformal flatness. To find the missing (spacetime-)orthogonal tangent vector to $\photo$, consider a curve $\mu(t)=(t,y(t))$ in $\photo$ with tangent vector $\dot{\mu}(t)=\partial_{t}+\dot{y}^{i}(t)\partial_{y^{i}}$. Let capital latin indices again run from $2,\dots,n$. Now decompose 
\begin{align}
\dot{y}(t)&=\rho(t) u(t,y(t))+\xi^{I}(t)e_{I}(t,y(t))\in\R^{n}.
\end{align}
By~\eqref{eq:spheres} and the fact that $\frac{d}{dt}\left\vert u(t,y(t))\right\vert_{\delta}^{2}=0$ for all $t\in T$, we find 
\begin{align}
\rho(t)&=\dot{y}(t)\cdot  u(t,y(t))=\dot{c}(t)\cdot u(t,y(t))+\dot{S}(t).
\end{align}
Hence, a future pointing tangent vector to $\photo$ orthogonal in the spacetime to all $X_{I}$ is given by
\begin{align}\label{eq:X1sphere}
X_{1}(t,y)&\definedas \partial_{t}+\left(\dot{c}(t)\cdot u(t,y)+\dot{S}(t)\right) u^{i}(t,y)\,\partial_{y^{i}},
\end{align}
so that we have constructed a smooth orthogonal tangent frame $\lbrace{X_{i}\rbrace}_{i=1}^{n}$ for $\photo$. Hence we can compute the outward (spacetime) unit normal to $\photo$ to be
\begin{align}\label{eq:etasphere}
\eta(t,y)&=\frac{\frac{\psi(s)}{\widetilde{N}(s)}\left(\dot{c}(t)\cdot u(t,y)+\dot{S}(t)\right)\partial_{t}+\frac{\widetilde{N}(s)}{\psi(s)} u^{i}(t,y)\partial_{y^{i}}}{\sqrt{\widetilde{N}^{2}(s)-\psi^{2}(s)\left(\dot{c}(t)\cdot u(t,y)+\dot{S}(t)\right)^{2}}}.
\end{align}
In other words, using that $\nu(t,y)=\psi^{-1}(s) u^{i}(t,y)\partial_{y^{i}}$, we have
\begin{align}\label{eq:alphaetasphere}
\alpha(t,y)&=\frac{\widetilde{N}(s)}{\sqrt{\widetilde{N}^{2}(s)-\psi^{2}(s)\left(\dot{c}(t)\cdot u(t,y)+\dot{S}(t)\right)^{2}}},\\\label{eq:betaetasphere}
\beta(t,y)&=\frac{\psi(s)\left(\dot{c}(t)\cdot u(t,y)+\dot{S}(t)\right)}{\widetilde{N}(s)\sqrt{\widetilde{N}^{2}(s)-\psi^{2}(s)\left(\dot{c}(t)\cdot u(t,y)+\dot{S}(t)\right)^{2}}}.
\end{align}
Extending the fields $e_{J}$ trivially in the radial direction, let us first collect the following explicit formulas arising from differentiating $e_{J}$ in direction $u$
\begin{align}\label{eq:trick1}
e_{I}^{i}(t,y)\left(e_{J},_{y^{i}}(t,y)\right)\cdot u(t,y)&=-\frac{e_{I}^{i}(t,y)\,e_{J}(t,y)\cdot y_{,y^{i}}}{S(t)}=-\frac{\delta_{IJ}}{S(t)},\\\label{eq:trick2}
u^{i}(t,y)\left(e_{J},_{y^{i}}(t,y)\right)\cdot u(t,y)&=-\frac{u^{i}(t,y)\,e_{J}(t,y)\cdot y_{,y^{i}}}{S(t)}=0.
\end{align}

Now, let us compute the second fundamental forms explicitly. Using~\eqref{eq:trick2}, we find that the umbilicity condition $\mathfrak{h}(X_{1},X_{J})=0$ gives
\begin{align*}
0&=u(t,y)\cdot\dot{e}_{J}(t,y)+\left(\dot{c}(t)\cdot u(t,y)+\dot{S}(t)\right) u_{i}(t,y) u^{k}(t,y)\\
&\quad\quad\times\left[e^{i}_{J\,,y^{k}}(t,y)+\frac{\psi'(s)}{s\psi(s)}\left(y_{k}\,e_{J}^{i}(t,y)+y\cdot e_{J}(t,y)\,\delta^{i}_{k}-y^{i}(e_{J})_{k}(t,y)\right)\right]\\
&\quad-\left(\dot{c}(t)\cdot u(t,y)+\dot{S}(t)\right)\frac{\widetilde{N}'(s)}{s\widetilde{N}(s)}\,y\cdot e_{J}(t,y)\\
&=\frac{\dot{c}(t)}{S(t)}\cdot e_{J}(t,y)+\left(\dot{c}(t)\cdot u(t,y)+\dot{S}(t)\right)\frac{1}{s}\left\{\frac{\psi'(s)}{\psi(s)}-\frac{\widetilde{N}'(s)}{\widetilde{N}(s)}\right\}\,y\cdot e_{J}(t,y)
\end{align*}
for all $J=2,\dots,n$. As $\lbrace{ u(t,y),e_{J}(t,y)\rbrace}_{J=2}^{n}$ is a $\delta$-orthonormal frame, this turns out to be equivalent to
\begin{align}\label{eq:mixedhspherical}
\begin{split}
-\frac{\dot{c}(t)}{S(t)}&=\left(\dot{c}(t)\cdot u(t,y)+\dot{S}(t)\right)\frac{1}{s}\left\{\frac{\psi'(s)}{\psi(s)}-\frac{\widetilde{N}'(s)}{\widetilde{N}(s)}\right\}y\\
&\quad-\left(\left(\dot{c}(t)\cdot u(t,y)+\dot{S}(t)\right)\frac{1}{s}\left\{\frac{\psi'(s)}{\psi(s)}-\frac{\widetilde{N}'(s)}{\widetilde{N}(s)}\right\}y\cdot  u(t,y)\right.\\
&\quad\quad\quad\left.+\,\frac{\dot{c}(t)\cdot u(t,y)}{S(t)}\right) u(t,y).
\end{split}
\end{align}
As $\surf(t)$ has dimension $n-1$ for all $t\in T$,~\eqref{eq:mixedhspherical} tells us by linear dependence considerations that $\dot{c}(t)=0$ for all $t\in T$. Consequentially,~\eqref{eq:mixedhspherical} simplifies to
\begin{align}\label{eq:mixedhsphericalsimple}
0&=\dot{S}(t)\left\{\frac{\psi'(s)}{\psi(s)}-\frac{\widetilde{N}'(s)}{\widetilde{N}(s)}\right\}\left(y
-(y\cdot u(t,y)) u(t,y)\right).
\end{align}
Assuming~\eqref{eq:unless}, the term in braces $\{\dots\}$ does not vanish and~\eqref{eq:mixedhsphericalsimple} implies that, for a fixed $t\in T$, $\dot{S}(t)=0$ or $y=(y\cdot u(t,y)) u(t,y)$ for all $y\in\surf(t)$. But $y=(y\cdot u(t,y)) u(t,y)$ for all $y\in\surf(t)$ is equivalent to $c=0$ and $S(t)=s$ for all $y\in\surf(t)$, again by linear dependence considerations and as $\surf(t)$ has dimension $n-1$. In other words, assuming~\eqref{eq:unless}, we now know that $c=0$ unless $S$ is constant along $\photo$ in which case a constant center $c\neq0$ is potentially possible.

Let us now continue with our computation of the conformal factor $\lambda$, using the simplification $\dot{c}(t)=0$ for all $t\in T$. We first treat the \textbf{case $\dot{S}(t)=0$ for all $t\in T$:} We know $S(t)=S=s$. Moreover,~\eqref{eq:alphaetasphere} and~\eqref{eq:betaetasphere} give $\alpha=1$, $\beta=0$ and by~\eqref{eq:etasphere} and~\eqref{eq:X1sphere} $\eta=\nu$ and $X_{1}=\partial_{t}$. Moreover, $u(t,y)=u(y)$ and $e_{J}(t,y)$ and hence $X_{J}$ are independent of $t$. By~\eqref{eq:trick1} and $y\cdot u(y)=c\cdot u(y)+S$ via~\eqref{eq:defu}, we find
\begin{align}
h(X_{I},X_{J})&=\frac{\psi(S)}{S}\left(1+\frac{\psi'(S)}{\psi(S)}\left(c\cdot u(y)+S\right)\right)\delta_{IJ}
\end{align}
so that, by~\eqref{eq:humbilic}, the photon surface umbilicity factor $\lambda$ satisfies
\begin{align}\label{eq:lambdahsphere}
\begin{split}
\lambda(t,y)&=\frac{1}{S\psi(S)}\left(1+\frac{\psi'(S)}{\psi(S)}\left(c\cdot u(y)+S\right)\right)
\end{split}
\end{align}
and thus in particular $\lambda(t,y)=\lambda(y)$ independent of $t$. Similarly, from $\mathfrak{h}(X_{1},X_{1})=\lambda p(X_{1},X_{1})$, we find
\begin{align}
\lambda(y)&=\frac{\widetilde{N}'(S)}{S\widetilde{N}(S)\psi(S)}\left(c\cdot u(y)+S\right)
\end{align}
and hence
\begin{align}
1+\left\{\frac{\psi'(S)}{\psi(S)}-\frac{\widetilde{N}'(S)}{\widetilde{N}(S)}\right\}\left(c\cdot u(y)+S\right)&=0
\end{align}
for all $y\in\surf(t)$ and all $t\in T$. As $\surf(t)$ has dimension $n-1$  and $S$ is constant, we conclude that $c\cdot u(y)$ is constant and hence by~\eqref{eq:defu} that $c\cdot y$ must be constant along $\photo$. Using again that $\surf(t)$ has dimension $n-1$, this leads to $c=0$ as desired. Hence, the isotropic radii $s=S_{*}$ for which this photon sphere can occur are the solutions of the implicit photon sphere equation
\begin{align}\label{eq:photonsphereisotropic}
1+\left\{\frac{\psi'(S_{*})}{\psi(S_{*})}-\frac{\widetilde{N}'(S_{*})}{\widetilde{N}(S_{*})}\right\}\,S_{*}&=0
\end{align}
provided such solutions exist.

Let us now treat the other \textbf{case $c=0$:} We find $X_{1}=\partial_{t}+\dot{S}(t)u^{i}(t,y)\partial_{y^{i}}$ by~\eqref{eq:X1sphere} and indeed $e_{J}(t,y)=e_{J}(y)$ is independent of $t$ as $u(t,y)=\frac{y}{S(t)}$. Moreover, $s=S(t)$ holds for all $(t,y)\in\photo$. Thus using~\eqref{eq:trick1}, we can compute
\begin{align}
h(X_{I},X_{J})&=\frac{\psi(S(t))}{S(t)}\left(1+\frac{\psi'(S(t))}{\psi(S(t))}S(t)\right)\delta_{IJ}
\end{align}
so that, by~\eqref{eq:humbilic}, the photon surface umbilicity factor $\lambda$ satisfies
\begin{align}\label{eq:lambdahspherecase2a}
\lambda(t,y)&=\frac{\widetilde{N}(S(t))\left(1+\frac{\psi'(S(t))}{\psi(S(t))}\,S(t)\right)}{S(t)\psi(S(t))\sqrt{\widetilde{N}^{2}(S(t))-\psi^{2}(S(t))\,\dot{S}^{2}(t)}}
\end{align}
from which we see that $\lambda(t,y)=\lambda(t)$ only depends on $t$. From the remaining umbilicity condition $\mathfrak{h}(X_{1},X_{1})=\lambda p(X_{1},X_{1})$, we obtain
\begin{align*}
\lambda(t)&=\frac{\widetilde{N}(S(t))\psi(S(t))}{\sqrt{\widetilde{N}^{2}(S(t))-\psi^{2}(S(t))\,\dot{S}^{2}(t)}^{\,3}}\\
&\quad\times \left(\frac{\widetilde{N}'(S(t))\widetilde{N}(S(t))}{\psi^{2}(S(t))}+\ddot{S}(t)+\left\{\frac{\psi'(S(t))}{\psi(S(t))}-\frac{2\widetilde{N}'(S(t))}{\widetilde{N}(S(t))}\right\}\dot{S}^{2}(t)\right)
\end{align*}
and can conclude that the implicit equation
\begin{align}\label{eq:photonsurfaceisotropic}
\begin{split}
&\left(1+\frac{\psi'(S(t))}{\psi(S(t))}\,S(t)\right)\left({\widetilde{N}^{2}(S(t))-\psi^{2}(S(t))\,\dot{S}^{2}(t)}\right)\\
&\quad=S(t)\widetilde{N}'(S(t))\widetilde{N}(S(t))\\
&\quad\quad+{S(t)\psi^{2}(S(t))}\left(\ddot{S}(t)+\left\{\frac{\psi'(S(t))}{\psi(S(t))}-\frac{2\widetilde{N}'(S(t))}{\widetilde{N}(S(t))}\right\}\dot{S}^{2}(t)\right)
\end{split}
\end{align}
holds for the isotropic radius profile $s=S(t)$.
\qed

\medskip
Now that we know that ``most'' photon surfaces in ``most'' spacetimes of class~$\mathcal{S}$ are spherically symmetric, let us gain a complementary perspective on these by relating spherically symmetric photon surfaces to null geodesics. The following notion of maximality will be useful for this endeavor.
\begin{Def}[Maximal photon surface]
Let $(\R\times \slice,\mathfrak{g})\in\mathcal{S}$ and let $\photo\hookrightarrow(\R\times \slice,\mathfrak{g})$ be a connected, spherically symmetric photon surface. We say that $P^{n}$ is \emph{maximal} if $P^{n}$ does not lie inside a strictly larger connected, spherically symmetric photon surface.
\end{Def}

If $\zeta\colon J\to \R\times \slice$ is a null geodesic in $(\R\times \slice,\mathfrak{g})\in\mathcal{S}$ defined on some interval $J\subseteq\R$ then the \emph{energy}  $E\definedas -\mathfrak{g}(\partial_{t},\dot{\zeta})$ is constant along $\zeta$. Exploiting the spherical symmetry of the spacetime, one can locally choose a linearly independent system of Killing vector fields $\{X_{K}\}_{K=1}^{n-1}$ representing the spherical symmetry and set $\ell_{K}\definedas \mathfrak{g}(X_{K},\dot{\zeta})$. Then $\ell\definedas \left(\Omega^{KL}\ell_{K}\ell_{L}\right)^{\frac{1}{2}}\geq0$ is called the \emph{(total) angular momentum of $\zeta$}. It is independent of the choice of the system $\{X_{K}\}_{K=1}^{n-1}$, and is constant along $\zeta$ as well. We will use the following definition of null geodesics generating a hypersurface.

\begin{Prop}\label{defprop:generating}
Let $(\R\times \slice,\mathfrak{g})\in\mathcal{S}$ and let $\zeta\colon J\to \R\times \slice$, $\zeta(s)=(t(s),r(s),\xi(s))$ be a (not necessarily maximal) null geodesic defined on some open interval $J\subseteq\R$ with angular momentum $\ell$, where $\xi(s)\in\mathbb{S}^{n-1}$. Then \emph{$\zeta$ generates the hypersurface $H^{n}_{\zeta}$} defined as
\begin{align*}
H^{n}_{\zeta}&\definedas\{(t,p)\in\R\times \slice\,\vert\,\exists s_{*}\in J, \xi_{*}\in\mathbb{S}^{n-1}:t=t(s_{*}),p=(r(s_{*}),\xi_{*})\},
\end{align*}
which is a smooth, connected, spherically symmetric hypersurface in $\R\times \slice$ which is timelike if $\ell>0$ and null if $\ell=0$.
\end{Prop}
\proof
The claim that $H^{n}_{\zeta}$ is a smooth, spherically symmetric hypersurface is verified by recalling that $r>0$ for all $(t,p)\in\R\times \slice$ so that $\zeta_{*}\in\mathbb{S}^{n-1}$ is unique and by realising that $t=t(s)$ is injective along the null geodesic $\zeta$ which shows that $s_{*}$ is unique. It is connected because the geodesic $\zeta$ is defined on an interval.

To show the claims about its causal character, set $e_{0}=\dot{\zeta}$ and extend it to all of $H^{n}_{\zeta}$ by making it invariant under the spherical symmetries. Thus, $e_{0}$ is a null tangent vector field to $H^{n}_{\zeta}$ orthogonal to each time-slice $\{t(s) = \text{const.}\}$ of $H^{n}_{\zeta}$. Now if $\ell=0$, it is easy to see that $\dot{\xi}=0$ along $\zeta$ and hence $e_{0}\perp_{\mathfrak{g}} \mathbb{S}^{n-1}$ or in other words $\mathfrak{g}(e_{0},X)=0$ for any vector field $X$ tangent to $H^{n}_{\zeta}$ whence $H^{n}_{\zeta}$ is null. Otherwise, if $\ell>0$, we will have $\dot{\xi}\neq0$ everywhere along $\zeta$ and hence $e_{0}$ will not be $\mathfrak{g}$-orthogonal to $\mathbb{S}^{n-1}$ and hence $\mathfrak{g}$ will induce a Lorentzian metric on~$H^{n}_{\zeta}$.
\qed

\medskip
With these concepts at hand, we can write down a characterization of spherically symmetric photon surfaces via generating null geodesics.
\begin{Prop}\label{prop:generating}
Let $(\R\times \slice,\mathfrak{g})\in\mathcal{S}$ and let $H^{n}\hookrightarrow(\R\times \slice,\mathfrak{g})$ be a connected, spherically symmetric, timelike hypersurface. Then $H^{n}$ is generated by a null geodesic $\zeta\colon J\to \R\times \slice$ if and only if $H^{n}$ is a photon surface. Moreover, $H^{n}$ is a maximal photon surface if and only if any null geodesic $\zeta\colon J\to \R\times \slice$ generating $H^{n}$ is maximal. 

The umbilicity factor $\lambda$ of a photon surface $P^{n}$ is related to the energy $E$ and angular momentum $\ell$ of its generating null geodesics by $\lambda=\frac{E}{\ell}$.
\end{Prop}
\proof
First, assume that $H^{n}$ is generated by a null geodesic $\zeta\colon J\to \R\times \slice$ so that $H^{n}=H^{n}_{\zeta}$ and observe that $H^{n}$ must then actually be ruled by the null geodesics arising by rotating $\zeta$ around $\mathbb{S}^{n-1}$. As $H^{n}$ is timelike by assumption, we know from Definition and Proposition \ref{defprop:generating} that the angular momentum $\ell$ of $\zeta$ satisfies $\ell>0$. Proceeding as above, let $\lbrace{e_{I}\rbrace}_{I=1}^{n-1}$ be a local orthonormal system tangent to $\mathbb{S}^{n-1}$ along $\zeta$, and set $e_{0}\definedas\dot{\zeta}$ and extend it to all of $H^{n}_{\zeta}$ by making it invariant under the spherical symmetries, so that $\lbrace{e_{0},e_{I}\rbrace}_{I=1}^{n-1}$ is a local frame along $\zeta$. Writing $\zeta(s)=(t(s),r(s),\xi(s))$ for $s\in J$ as before, we can write 
\begin{align}
\dot{\zeta}(s)&=\dot{t}(s)\partial_{t}+\dot{r}(s)\partial_{r}+\dot{\xi}(s)
\end{align}
for $s\in J$ so that $\zeta$ being a null curve is equivalent to 
\begin{align}
-f(r)\dot{t}^{2}+\frac{\dot{r}^{2}}{f(r)}+r^{2}\vert\dot{\xi}\vert^{2}_{\Omega}&=0
\end{align}
along $\zeta$. From this and the definition of energy $E$ and angular momentum $\ell$, we obtain
\begin{align}\label{eq:Edt}
\dot{t}&=\frac{E}{f(r)},\\\label{eq:rdr}
\dot{r}^{2}&=E^{2}-\frac{\ell^{2}f(r)}{r^{2}},
\end{align}
\begin{align}
\vert\dot{\xi}\vert^{2}_{\Omega}&=\frac{\ell^{2}}{r^{4}}
\end{align}
along $\zeta$. To compute the second fundamental form $\mathfrak{h}$ of $H^{n}\hookrightarrow(\R\times\slice,\mathfrak{g})$, let us first compute the outward (growing $r$) spacelike unit normal $\eta$ to $H^{n}$. By spherical symmetry, $\eta$ must be a linear combination of $\partial_{t}$ and $\partial_{r}$, with no angular contribution. It must also be orthogonal to $e_{0}=\dot{\zeta}$. From this, we find
\begin{align}
\eta&=\frac{1}{\ell}\left(\frac{r\dot{r}}{f(r)}\partial_{t}+Er\partial_{r}\right),
\end{align}
where we have used \eqref{eq:Edt}, \eqref{eq:rdr}. For $I=1,\dots,n-1$, we find by a direct computation (exploiting spherical symmetry) that
\begin{align}
\nabla_{e_{I}}\eta&=\frac{E}{\ell}\,e_{I}
\end{align}
and hence
\begin{align}
\mathfrak{h}(e_{I},e_{\beta})&=\mathfrak{g}(\nabla_{e_{I}}\eta,e_{\beta})=\frac{E}{\ell}\,\mathfrak{g}(e_{I},e_{\beta})
\end{align}
for $I=1,\dots,n-1$ and $\beta=0,\dots,n-1$. On the other hand, smoothly extending $e_0$ to a neighborhood of $H^n$, one finds
\begin{align}
\mathfrak{h}(e_{0},e_{0})&=\mathfrak{g}(\nabla_{e_{0}}\eta,e_{0})=-\mathfrak{g}(\nabla_{e_{0}}e_{0},\eta)=0
\end{align}
holds as $e_{0}=\dot{\zeta}$ and $\zeta$ is a geodesic. Hence by symmetry of $\mathfrak{h}$ and $\mathfrak{g}$
\begin{align}
\mathfrak{h}(e_{\alpha},e_{\beta})&=\frac{E}{\ell}\mathfrak{g}(e_{\alpha},e_{\beta})
\end{align}
for $\alpha,\beta=0,\dots,n-1$ so that $H^{n}$ is totally umbilic with umbilicity factor $\lambda=\frac{E}{\ell}$.

Conversely, assume that $H^{n}$ is a spherically symmetric photon surface. Let $q\in H^{n}$ and let $X\in T_{q}H^{n}$ be any null tangent vector. Let $\zeta\colon J_{\text{max}}\ni0\to \R\times \slice$ be the maximal null geodesic with $\zeta(0)=q$, $\dot{\zeta}(0)=X$. As $H^{n}$ is totally umbilic, $\zeta$ must remain tangent to $H^{n}$ on some maximal open interval $0\in J\subseteq J_{\text{max}}$ by Proposition~\ref{prop:umbilic}. By spherical symmetry, $H^{n}$ is in fact generated by $\zeta\vert_{J}$. 

It remains to discuss the claim about maximality which is a straightforward consequence of the above argument about the choice of the interval $J\subseteq J_{\text{max}}$.\qed

\begin{Rem}
In view of effective one-body dynamics (see for example~\cite{BD}), it may be of interest to point out that a spherically symmetric photon surface $P^{n}$ in a static, spherically symmetric spacetime with metric of the form
\begin{align*}
\mathfrak{g}&=-Gdt^{2}+\frac{1}{f}dr^{2}+r^{2}\,\Omega,
\end{align*}
where $G=G(r)$, $f=f(r)$, are smooth, positive functions on an open interval $I\subseteq(0,\infty)$, will have constant umbilicity factor $\lambda$ if and only if $G=\kappa f$ along $P^{n}$ for some $\kappa>0$, with $\lambda=\frac{Ef(r)}{\ell \,G(r)}$ characterizing the umbilicity factor along any generating null geodesic with energy $E$ and angular momentum $\ell$.
\end{Rem}

\begin{Rem}
Proposition \ref{prop:generating} allows us to conclude existence of (maximal) spherically symmetric photon surfaces in spacetimes of class $\mathcal{S}$ from existence of (maximal) null geodesics.
\end{Rem}

A different view on spherically symmetric photon surfaces in spacetimes of class $\mathcal{S}$ can be gained by lifting them to the phase space (i.e. to the cotangent bundle). We will end this section by proving a foliation property of the null section of phase space by maximal photon surfaces and by so-called maximal principal null hypersurfaces.

\begin{Def}[(Maximal) Principal null hypersurfaces]\label{def:principal}
Let $(\R\times \slice,\mathfrak{g})\in\mathcal{S}$, $\zeta\colon J\to\R\times\slice$ a null geodesic, and $\ell$ its angular momentum. Then $\zeta$ is called a \emph{principal null geodesic} if $\ell=0$. A hypersurface $H^{n}_{\zeta}$ generated by a (maximal) principal null geodesic $\zeta$ will be called a \emph{(maximal) principal null hypersurface}.
\end{Def}

Recall that from Definition and Proposition \ref{defprop:generating}, we know that principal null hypersurfaces are spherically symmetric and indeed null. Arguing as in the proof of Proposition \ref{prop:generating}, one sees that a maximal principal null hypersurface will be maximal in the sense that it is not contained in any strictly larger principal null hypersurface. In particular, if one generating null geodesic is maximal, then all of them are. Finally, principal null hypersurfaces are connected by definition.

With these considerations at hand, let us prove the following foliation property of the null section of the phase space of any spacetime of class $\mathcal{S}$. To express this, we will canonically lift the null bundles over the involved spherically symmetric photon surfaces and principal null hypersurfaces to the null section of the phase space. Recall that maximal spherically symmetric photon surfaces are connected by definition.

\begin{Prop}
Let $(\R\times \slice,\mathfrak{g})\in\mathcal{S}$. Then the null section of the phase space of $(\R\times\slice,\mathfrak{g})$, $\mathcal{N}\definedas\{\omega\in T^{*}(\R\times\slice)\,\vert\, \mathfrak{g}(\omega^{\#},\omega^{\#})=0\}$, is foliated by the canonical lifts $\mathcal{N}(P^{n})\subset\mathcal{N}$ of the null bundles over all maximal spherically symmetric photon surfaces $P^{n}$ and of the canonical lifts $\mathcal{N}(H^{n})\subset\mathcal{N}$ of the null bundles over all maximal principal null hypersurfaces $H^{n}$.
\end{Prop}
\proof
Consider $\omega\in\mathcal{N}$ and let $\zeta\colon J\ni0\to\R\times\slice$ be the unique maximal null geodesic satisfying  $\dot{\zeta}(0)=\omega^{\#}$. Let $H^{n}_{\zeta}$ be the hypersurface of $\R\times\slice$ generated by $\zeta$ and note that $\omega\in\mathcal{N}(H^{n}_{\zeta})$ holds by construction. By Definition and Proposition \ref{defprop:generating}, Proposition \ref{prop:generating}, and Definition \ref{def:principal}, we know that $H^{n}_{\zeta}$ is a maximal photon surface if the angular momentum $\ell$ of $\zeta$ is positive, and a maximal principal null hypersurface if $\ell=0$ vanishes. As $\ell\geq0$, this shows that $\omega$ is either contained in at least one canonical lift of a maximal photon surface or contained in at least one canonical lift of a maximal principal null hypersurface. Furthermore, $\omega$ cannot lie in the canonical lifts of two different maximal spherically symmetric photon surfaces because these would both be generated by $\zeta$ and hence coincide by Proposition~\ref{prop:generating}. Finally, $\omega$ cannot lie in the canonical lifts of two different maximal principal null hypersurfaces which can be seen by repeating the arguments in the proof of Proposition~\ref{prop:generating}.
\qed

 \section{A rigidity result for photon surfaces with equi\-potential time-slices}\label{sec:equi}\label{sec:unique}
As discussed in the previous section, the \schild spacetime of mass $m>0$ in $n+1$ dimensions possesses not only the well-known photon sphere at $r=(nm)^{\frac{1}{n-2}}$ but also many other photon surfaces. Except for the planar ones, all of these \schild photon surfaces are spherically symmetric and thus in particular equipotential as defined in Section~\ref{sec:prelim}. In this section, we will prove the following theorem which can be considered complementary to Corollary~\ref{coro:schwarzphoto} in the context of static, vacuum, asymptotically flat spacetimes.

\begin{Thm}\label{thm:main}
Let $(\mathfrak{L}^{n+1},\mathfrak{g})$ be a static, vacuum, asymptotically isotropic spacetime of mass $m$. Assume that $(\mathfrak{L}^{n+1},\mathfrak{g})$ is geodesically complete up to its inner boundary $\partial\mathfrak{L}$, which is assumed to be a (possibly disconnected) photon surface, $\partial\mathfrak{L}\asdefined\photo$. Assume in addition that $\photo$ is equipotential,  outward directed, and has compact time-slices $\surf(t)=\photo\cap\slice(t)$. Then $(\mathfrak{L}^{n+1},\mathfrak{g})$ is isometric to a suitable piece of the \schild spacetime of mass $m$, and in fact $m>0$. In particular, $\photo$ is connected, and is (necessarily) a spherically symmetric photon surface in \schild spacetime. 
\end{Thm}

The proof relies on the following theorem by the first named author.
\begin{Thm}[{\!\!\cite{ndimunique}}]\label{thm:hidim}
Assume $n\geq3$ and let $\slice$ be a smooth, connected, $n$-dimensional manifold with non-empty, possibly disconnected, smooth, compact inner boundary $\partial M=\dt\cup_{i=1}^{I}\surf_{i}$. Let $g$ be a smooth Riemannian metric on $\slice$. Assume that $(\slice,g)$ has non-negative scalar curvature
\begin{align*}
\Scal\geq0,
\end{align*} and that it is geodesically complete up to its inner boundary $\partial M$. Assume in addition that $(\slice,g)$ is asymptotically isotropic with one end of mass $m\in\R$. Assume that the inner boundary $\partial M$ is umbilic in $(\slice,g)$, and that each component $\surf_{i}$ has constant mean curvature $H_{i}$ with respect to the outward pointing unit normal~$\nu_{i}$. Assume furthermore that there exists a function $u\colon \slice\to\R$ with $u>0$ away from $\partial M$ which is smooth and harmonic on $(\slice,g)$, 
\begin{align*}
\triangle u=0.
\end{align*} We ask that $u$ is asymptotically isotropic of the same mass $m$, and such that $u\vert_{\Sigma_{i}^{n-1}}\equiv:u_{i}$ and the normal derivative of $u$ across $\surf_i$, $\nu_{i}(u)\vert_{\Sigma_{i}^{n-1}}\equiv:\nu(u)_{i}$, are constant on each $\surf_{i}$. Finally, we assume that for each $i=1,\dots,I$,  we are either in the \emph{semi-static horizon case}
\begin{align}\label{eq:H0}
H_{i}=0,\quad u_i=0, \quad \nu(u)_{i}\neq0,
\end{align}
or in the \emph{true CMC case} $H_{i}>0$, $u_{i}>0$, and there exists $c_{i}>\frac{n-2}{n-1}$ such that
\begin{align}\label{eq:scali}
\Scal_{\sigma_{i}}&=c_i H_{i}^{2},\\\label{eq:nuui}
2\nu(u)_{i}&=\left(c_{i}-\frac{n-2}{n-1}\right)H_{i}u_{i},
\end{align}
where $\Scal_{\sigma_{i}}$ denotes the scalar curvature of $\surf_{i}$ with respect to its induced metric $\sigma_{i}$.

Then $m>0$ and $(\slice,g)$ is isometric to a suitable piece $(\widetilde{M}^{n}_{m}\setminus B_{S}(0),\widetilde{g}_{m})$ of the (isotropic) \schild manifold of mass $m$ with $S\geq s_{m}$. Moreover, $u$ coincides with the restriction of $\widetilde{u}_{m}$ (up to the isometry) and the isometry is smooth.
\end{Thm}

\begin{Rem}[Generalization]\label{rem:gen}
Our proof of Theorem \ref{thm:main} makes use of the static vacuum Einstein equations,~\eqref{SMEvac1},~\eqref{SMEvac2} and~\eqref{SMEvac3}.
In fact, as we will see in the proof below, it is sufficient to ask that the vacuum Einstein equations hold in a neighborhood of $\photo$; outside this neighborhood it suffices that $\triangle N=0$ and that the dominant energy condition $R\geq0$ holds. 
\end{Rem}

\begin{Rem}[Multiple ends]\label{rem:multiple}
As Theorem \ref{thm:hidim} generalizes to multiple ends (see~\cite{ndimunique}), Theorem \ref{thm:main} also readily applies in the case of multiple ends satisfying the decay conditions~\eqref{eq:AI} and~\eqref{eq:NAI} with potentially different masses $m_{i}$ in each end $E_{i}$. Note that, in each end $E_{i}$, it is necessary that both $g_{ij}$ and $N$ have the same mass $m_{i}$ in their expansions.
\end{Rem}

\begin{Rem}[Discussion of $\eta(N) > 0$] The assumption that $\photo$ is outward directed, $\eta(N) > 0$ (hence $dN \ne 0$) along  $\photo$,  can be removed if instead, one assumes that $m>0$ a priori and that $\photo$ is connected. Using the Laplace equation $\triangle N=0$ and the divergence theorem as well as the asymptotics~\eqref{eq:AI} and~\eqref{eq:NAI}, one computes
\begin{align*}
\frac{1}{\omega_{n-1}}\int_{\surf}\nu(N)\,dA&=m>0,
\end{align*}
where $\surf\definedas \photo\cap\lbrace t=\text{const}\rbrace$ and $\omega_{n-1}$ is the volume of $(\mathbb{S}^{n-1},\Omega)$, see \cite[Definition~4.2.1]{CDiss} for the $n=3$ case; the argument is identical in higher dimensions. From this and connectedness of $\surf$, we can deduce that $\nu(N)>0$ and thus in particular $\eta(N) > 0$ at least in an open subset of $\surf$. However, we will see in the proof of Theorem \ref{thm:main} that this necessarily implies that $\nu(N)\equiv\text{const}>0$ in this open neighborhood (noting that all computations performed there are purely local). As $\surf$ is connected, we obtain $\nu(N)\equiv\text{const}>0$ and thus in particular $\eta(N) > 0$ everywhere on $\surf$, see Equation~\eqref{eq:eta} below.
\end{Rem}

\proof
 We write $\slice(t)$ for the time-slice $\{t\} \times M$ (cf.  Remark \ref{rem:slice}), and consider each connected component $\photo_{i}$, $i=1,\dots I$ of $\photo$ separately. For the component of $\photo_{i}$ under consideration, let $\surf_{i}(t)\definedas \photo_{i}\cap \slice(t)$. We will drop the explicit reference to $i$ in what follows and only start reusing it toward the end of the proof, where we bring in global arguments.

Let $\nu$ denote the outward unit normal to $\surf(t)\hookrightarrow (\slice(t),g)$, pointing to the asymptotically isotropic end. Let $\eta$ denote the outward unit normal of $\photo\hookrightarrow (\mathbb{R}\times \slice,\mathfrak{g})$. As $N\asdefined u(t)$ on $\surf(t)$ and because we assumed $\eta(N) > 0$, and hence 
$\nu(N) >0$ (see~\eqref{eq:eta} below) on $\photo$, we have
\begin{align}
\nu=&\frac{\grad N}{\vert \grad N\vert}=\frac{\grad N}{\vert dN\vert}=\frac{\grad N}{\nu(N)}.
\end{align}
Now let $\mu(s)=(s,x(s))$ be a curve in $\photo$, i.e.\ $N\circ\mu(s)=u(s)$. This implies by chain rule that $dN(\dot{\mu})=\dot{u}$. If $\dot{\mu}(t)\perp \surf(t)$, the tangent vector of $\mu$ can be computed explicitly as 
\begin{align}\label{eq:Z}
Z\definedas\dot{\mu}=\d_{t}+\dot{x}=\d_{t}+\frac{\dot{u}}{\nu(N)}\nu\in\Gamma(T\photo).
\end{align}
Expressed in words, $Z$ is the vector field going ``straight up'' along $\photo$.\\

The explicit formula~\eqref{eq:Z} for $Z$ allows us to explicitly compute the spacetime unit normal $\eta$ to $\photo$, too: It has to be perpendicular to $Z$ and its projection onto $\slice$ has to be proportional to $\nu$. From this, we find
\begin{align}\label{eq:eta}
\eta=\frac{\nu+\frac{\dot{u}}{u^{2}\nu(N)^{2}}\d_{t}}{\sqrt{1-\frac{\dot{u}^{2}}{u^{2}\nu(N)^{2}}}}.
\end{align}
From umbilicity of $\photo\hookrightarrow (\mathbb{R}\times \slice,\mathfrak{g})$, it follows that the corresponding second fundamental form $\mathfrak{h}$ of $\photo\hookrightarrow (\mathbb{R}\times \slice,\mathfrak{g})$ satisfies
\begin{align}\label{eq:hfrak}
\mathfrak{h}=\frac{1}{n}\mathfrak{H}\, p,
\end{align}
where $p$ is the induced metric on $\photo$ and $\mathfrak{H}\definedas\tr_{p}\mathfrak{h}$. From Proposition 3.3 in~\cite{CederPhoto}, we know that $\mathfrak{H}\equiv\text{const}$. Equation~\eqref{eq:hfrak} implies in particular, that, for any tangent vector fields $X,Y\in\Gamma(T\surf(t))$, we have
\begin{align}\label{eq:umbilic}
\mathfrak{h}(X,Y)&=\frac{1}{n}\,\mathfrak{H}\,\sigma(X,Y),
\end{align}
where $\sigma$ denotes the induced metric on $\surf(t)$. Now extend $X,Y$ arbitrarily smoothly along $\photo$ such that they remain tangent to $\surf_{t}$. We compute
\begin{align}\nonumber
\mathfrak{h}(X,Y)&=-\mathfrak{g}\left(^{\mathfrak{g}}\nabla_{X}Y,\eta\right)\\\nonumber
&=-\frac{1}{\sqrt{1-\frac{\dot{u}^{2}}{u^{2}\nu(N)^{2}}}}\,\mathfrak{g}\left(^{\mathfrak{g}}\nabla_{X}Y,\nu+\frac{\dot{u}}{u^{2}\nu(N)^{2}}\d_{t}\right)\\\nonumber
&=-\frac{1}{\sqrt{1-\frac{\dot{u}^{2}}{u^{2}\nu(N)^{2}}}}\lbrace{\mathfrak{g}\left(^{\mathfrak{g}}\nabla_{X}Y,\nu\right)+\frac{\dot{u}}{u^{2}\nu(N)^{2}}\,\mathfrak{g}\left(^{\mathfrak{g}}\nabla_{X}Y,\d_{t}\right)\rbrace}\\\nonumber
&=-\frac{1}{\sqrt{1-\frac{\dot{u}^{2}}{u^{2}\nu(N)^{2}}}}\lbrace{\mathfrak{g}\left(^{\mathfrak{g}}\nabla_{X}Y,\nu\right)-\frac{\dot{u}}{u^{2}\nu(N)^{2}}{K\left(X,Y\right)\rbrace}}\\\nonumber
&\stackrel{K=0}{=}-\frac{1}{\sqrt{1-\frac{\dot{u}^{2}}{u^{2}\nu(N)^{2}}}}\,\mathfrak{g}\left(^{\mathfrak{g}}\nabla_{X}Y,\nu\right)\\\nonumber
&\stackrel{\mathfrak{g}\text{ static}}{=}-\frac{1}{\sqrt{1-\frac{\dot{u}^{2}}{u^{2}\nu(N)^{2}}}}\,g\left(^g\nabla_{X}Y,\nu\right)\\\label{eq:h}
&=\frac{1}{\sqrt{1-\frac{\dot{u}^{2}}{u^{2}\nu(N)^{2}}}}\,h(X,Y),
\end{align}
where $K=0$ denotes the second fundamental form of $\slice(t)\hookrightarrow (\mathbb{R}\times \slice,\mathfrak{g})$ and $h$ denotes the second fundamental form of $\surf(t)\hookrightarrow (\lbrace{t}\rbrace\times \slice,g)$. In particular, $\surf(t)\hookrightarrow (\lbrace{t}\rbrace\times \slice,g)$ is umbilic and its mean curvature $H$ inside $\slice$ can be computed as
\begin{align}\label{eq:H}
H&=\frac{n-1}{n}\,\mathfrak{H}\,\sqrt{1-\frac{\dot{u}^{2}}{u^{2}\nu(N)^{2}}}
\end{align}
when combining~\eqref{eq:umbilic} with~\eqref{eq:h}. Of course then $h=\frac{1}{n-1}H\sigma$. We will now proceed to show that $H$ and $\nu(N)$ are constant for each fixed $t$. Consider first the contracted Codazzi equation for $\surf(t)\hookrightarrow (\lbrace{t}\rbrace\times \slice,g)$. It gives us
\begin{align}\label{eq:Codazzi}
\Ric_{g}(X,\nu)&=\frac{n-2}{n-1} \,X(H).
\end{align}
On the other hand, using the static equation gives us
\begin{align}\nonumber
X(\nu(N))&\stackrel{N\equiv u(t)}{=}X(\nu(N))-{}^{g}\nabla_{X}\nu (N)\\\nonumber
&\;\;\;={}^{g}\nabla^{2}N(X,\nu)\\\nonumber
&\;\;\;=N\,{}^{g}\!\Ric(X,\nu)\\\label{eq:XH}
&\stackrel{\eqref{eq:Codazzi}}{=}\frac{n-2}{n-1}\,u\,X(H).
\end{align}
Furthermore,~\eqref{eq:H} allows us to compute
\begin{align}\nonumber
X(H)&=\frac{(n-1)\,\mathfrak{H}}{n\,\sqrt{1-\frac{\dot{u}^{2}}{u^{2}\nu(N)^{2}}}}\,\frac{\dot{u}^{2}}{u^{2}\,\nu(N)^{3}}\, X(\nu(N))\\
&\stackrel{\eqref{eq:XH}}{=}\frac{(n-2)\,\dot{u}^{2}\,\mathfrak{H}}{nu\,\nu(N)^{3}\,\sqrt{1-\frac{\dot{u}^{2}}{u^{2}\nu(N)^{2}}}}\, X(H).
\end{align}
Assume $X(H)\neq0$ in some open subset $U\subset\surf(t)$. Then in $U$, we have
\begin{align*}
nu\,\nu(N)^{3}\,\sqrt{1-\frac{\dot{u}^{2}}{u^{2}\nu(N)^{2}}}&=(n-2)\,\dot{u}^{2}\,\mathfrak{H}\\
\Leftrightarrow\quad\quad\nu(N)^{6}\left(1-\frac{\dot{u}^{2}}{u^{2}\nu(N)^{2}}\right)&=\frac{(n-2)^{2}\,\dot{u}^{4}\,\mathfrak{H}^{2}}{n^{2}\,u^{2}}\\
\Leftrightarrow\quad\quad\nu(N)^{6}-\frac{\dot{u}^{2}}{u^{2}}\,\nu(N)^{4}-\frac{(n-2)^{2}\,\dot{u}^{4}\,\mathfrak{h}^{2}}{n^{2}\,u^{2}}&=0.
\end{align*}
This is a polynomial equation for $\nu(N)$ with coefficients that only depend on $t$. As a consequence, $\nu(N)$ has to be constant in $U$. However, from~\eqref{eq:XH}, we know that then also $H$ has to be constant in $U$, a contradiction to $X(H)\neq0$ in $U$. Thus, $H$ and by~\eqref{eq:XH} also $\nu(N)$ are constants along $\surf(t)$ and only depend on $t$. From now on, we will drop the explicit reference to $t$ and also go back to using $N$ instead of $u(t)$ as the remaining part of the proof applies to each $t$ separately.\\

Thus, each $\surf=\surf(t)$ is an umbilic, CMC, equipotential surface in $\slice$ with $\nu(N)$ constant, too. From the usual decomposition of the Laplacian on functions and the static vacuum equation~\eqref{SMEvac3}, we find
\begin{align*}
0&=\triangle N=\triangle_{\sigma} N+N\Ric(\nu,\nu)+H\nu(N)=N\Ric(\nu,\nu)+H\nu(N)
\end{align*}
so that $\Ric(\nu,\nu)=-\frac{H\nu(N)}{N}$ must also be constant along $\surf$.

Plugging this into the contracted Gau{\ss} equation and using $\Scal=0$, we obtain
\begin{align}\nonumber
-2\Ric(\nu,\nu)&=\Scal_{\sigma}-\frac{n-2}{n-1}H^{2}\\\label{eq:scalformula}
\Leftrightarrow\quad\quad \Scal_{\sigma}&=\frac{2H\nu(N)}{N}+\frac{n-2}{n-1}H^{2}.
\end{align}
This shows that $(\surf,\sigma)$ also has constant scalar curvature. Now define the constant $c>\frac{n-2}{n-1}$ by
\begin{align}\label{eq:c}
c\definedas\frac{n-2}{n-1}+\frac{2\nu(N)}{NH}.
\end{align}
Together with~\eqref{eq:scalformula}, this definition of $c$ ensures
\begin{align}\label{constraint1}
R_{\sigma}&=c H^{2},\\\label{constraint2}
2\nu(N)&=\(c-\frac{n-2}{n-1}\) HN.
\end{align}

Let us summarize: Fix $t$ and keep dropping the explicit reference to it. Then each component of $(\surf,\sigma)\hookrightarrow(\slice,g)$ is umbilic, CMC, equipotential, has constant scalar curvature and constant $\nu(N)$, and all these constants together satisfy Equations~\eqref{constraint1} and~\eqref{constraint2} with constant $c$ given by~\eqref{eq:c} which is potentially different for each component of $\surf$. Recall the assumption that $(\slice,g)$ is geodesically complete up to its inner boundary so in particular $(\slice\setminus K,g)$ is geodesically complete up to $\surf$, where $K$ is the compact set such that $\surf=\partial\(\slice\setminus K\)$. Moreover, $(\slice\setminus K,g)$ satisfies the static vacuum equations. Altogether, these facts ensure that 
Theorem~\ref{thm:hidim} applies. Thus $(\slice\setminus K,g)$ is isometric to a spherically symmetric piece of the spatial \schild manifold of mass $m$ given by the asymptotics~\eqref{eq:AI},~\eqref{eq:NAI}, and $N$ corresponds to the \schild lapse function of the same mass $m$ under this isometry. The area radius $r$ of the inner boundary $\surf$ in the spatial \schild manifold is determined by $R_{\sigma}=\frac{(n-1)(n-2)}{r^{2}}$.\footnote{We point out in connection with Remarks \ref{rem:multiple} and \ref{rem:gen} that 
the assumptions of Theorem \ref{thm:hidim} keep being met if we start with several ends and the static vacuum equations only holding near ${\Sigma}^{n-1}$, with $\triangle N=0$ everywhere in $\slice$ as all the above computations and arguments were purely local near $\photo$.}

Thus, recalling the dependence on $t$, the manifold $(\lbrace t\rbrace \times (\slice\setminus K(t)),g)$ outside the photon surface time-slice $\surf(t)$ is isometric to a piece of the spatial \schild manifold of mass $m$ with inner boundary area radius $r(t)$. In particular, $\surf(t)$ is a connected sphere and we find $m>0$ (as per Theorem 1.1 in~\cite{ndimunique}). Recombining the time-slices $\surf(t)$ to the photon surface $\photo$, this shows that the part of the spacetime $(\mathbb{R}\times \slice,\mathfrak{g})$ lying outside the photon surface $\photo$ is isometric to a piece of the \schild spacetime of mass $m$, and $m>0$ necessarily. Moreover, $\photo$ is connected and its isometric image in the \schild spacetime is spherically symmetric with radius profile~$r(t)$. Moreover, $\photo$ is connected and its isometric image in the \schild spacetime is spherically symmetric with radius profile~$r(t)$. \qed

\paragraph{Acknowledgements.} The authors would like to thank Gary Gibbons, Sophia Jahns, Volker Perlick, Volker Schlue, Olivia Vi\v{c}\'anek Mart\'inez, and Bernard Whiting for helpful comments and questions. 

The authors would like to extend thanks to the Mathematisches Forschungsinstitut Oberwolfach, the University of Vienna, and the Tsinghua Sanya International Mathematics Forum for allowing us to collaborate in stimulating environments.

The first named author is indebted to the Baden-W\"urttemberg Stiftung for the financial support of this research project by the Eliteprogramme for Postdocs. The work of the first named author is supported by the Institutional Strategy of the University of T\"ubingen (Deutsche Forschungsgemeinschaft, ZUK 63) and by the focus program on Geometry at Infinity (Deutsche Forschungsgemeinschaft,  SPP 2026).  The work of the second named author was partially supported by NSF grants  DMS-1313724 and DMS-1710808. 

\bibliographystyle{amsplain}
\bibliography{photon-surfaces}

\end{document}